\documentclass[letterpaper, 10 pt, conference]{ieeeconf}
\IEEEoverridecommandlockouts

\usepackage{graphicx}
\usepackage{cite}
\usepackage{authblk}
\usepackage{subfigure}
\usepackage{multirow}
\usepackage{color}
\usepackage{tikz} \usetikzlibrary{arrows,calc,decorations.pathreplacing}

\usepackage{amsmath,amssymb,amsfonts}
\usepackage{algorithmic}
\usepackage{graphicx}
\usepackage{textcomp}
\usepackage{xcolor}
\def\BibTeX{{\rm B\kern-.05em{\sc i\kern-.025em b}\kern-.08em
    T\kern-.1667em\lower.7ex\hbox{E}\kern-.125emX}}
\usepackage{url}
\usepackage{tikz}
\usetikzlibrary{arrows,calc,decorations.markings}
\usepackage[siunitx]{circuitikz}
\graphicspath{{Figures/}}
\usepackage[strings]{underscore} 
\usepackage{todonotes}

\newtheorem{theorem}{Theorem}

\newtheorem{proposition}{Proposition}

\title{An Introduction to Pursuit-evasion Differential Games}
\author{Isaac E. Weintraub, Meir Pachter, and Eloy Garcia
\thanks{I. Weintraub is with the Aerospace Systems Directorate, Air Force Research Laboratory, Wright-Patterson AFB, OH 45433. \ttfamily{isaac.weintraub.1@us.af.mil}}
\thanks{M. Pachter is with the Department of Electrical Engineering, Air Force Institute of Technology, Wright-Patterson AFB, OH 45433. \ttfamily{meir.pachter@afit.edu}}
\thanks{E. Garcia is with the Control Science Center of Excellence, Air Force Research Laboratory, Wright-Patterson AFB, OH 45433. \ttfamily{eloy.garcia.2@us.af.mil}}
}

\begin{document}
\maketitle 
\begin{abstract}
Pursuit and evasion conflicts represent challenging problems with important applications in aerospace and
      robotics. In pursuit-evasion problems, synthesis of
      intelligent actions must consider the adversary's potential
      strategies. Differential game theory provides an adequate
      framework to analyze possible outcomes of the conflict
      without assuming particular behaviors by the opponent. This
      article presents an organized introduction of
      pursuit-evasion differential games with an overview of
      recent advances in the area. First, a summary of the seminal work is outlined, highlighting important contributions. Next, more recent results are described by employing a classification based on the number of players: one-pursuer-one-evader, N-pursuers-one-evader, one-pursuer-M-evaders, and N-pursuer-M-evader games. In each scenario, a brief summary of the literature is presented. Finally, two representative
      pursuit-evasion differential games are studied in detail:
      the two-cutters and fugitive ship differential game and the
      active target defense differential game. These problems
      provide two important applications and, more importantly,
      they give great insight into the realization of cooperation
      between friendly agents in order to form a team and defeat
      the adversary. 
\end{abstract}

\section{Introduction} \label{sec:Intro}

Pursuit-evasion problems provide a general framework that mathematically formalizes important applications in different areas such as surveillance, navigation, analysis of biological behaviors, and conflict and combat operations. Pursuit-evasion sets up (in the simplest case) two players or autonomous agents against each other; generalizations are typical in the sense of multiple players divided into two teams -- the pursuer team against the evader team.  The main purpose is to provide strategies that enable an autonomous agent to perform a set of actions against the opponent, for instance, the pursuer aims at determining a strategy that will result in capture or interception of the evader. 

It is common to approach pursuit-evasion problems by imposing certain assumptions on the behavior of the opponent \cite{Sprinkle04,EarlDandrea07}. 
However, many pursuit-evasion scenarios must address the presence of an intelligent adversary which does not abide by a restricted set of actions. The desire to design strategies that optimize a certain criterion against the worst possible actions of the opponent and that also provide robustness with respect to all possible behaviors of the same opponent, led to the emergence of differential game theory \cite{Isaacs1965Differential}.
The central problem in pursuit-evasion differential games is the synthesis of saddle-point strategies that provide guaranteed performance for each team regardless of the actual strategies implemented by the adversary.

This paper serves as an introduction into the area of two-team pursuit-evasion differential games with a focus on multi-player problems, that is, games with three or more players where team cooperation is implicitly required. An outline of seminal work on the area of differential games which highlights important contributions and potential applications is given in Section \ref{sec:litSeminal}. This is followed by an overview of important contributions in the field using a suitable categorization based on number of players in Section \ref{sec:litPE}. Cooperation between agents of the same team is emphasized both in theoretical and practical terms. From a theoretical point of view, the saddle-point solution of a differential game with multiple players in one or both teams requires the strictest level of cooperation. From a practical perspective, cooperation improves team performance compared to the individual efforts. The paper culminates with two representative case studies in Sections \ref{sec:TwoCutters} and \ref{sec:litATDDG} which illustrate both the tools available to synthesize and verify optimal strategies and the important cooperation aspect within multi-player differential games.

\section{Seminal Work} 
\label{sec:litSeminal}
The development of differential games started with the works of \cite{Isaacs1951Games, Isaacs1954DG1, Isaacs1954DG2, Isaacs1954DG3, Isaacs1955DG4, Isaacs1965Differential}. In these publications, Isaacs outlined the idea of posing problems in a dynamic game-theoretic framework; he called this paradigm ``Differential Games''. In his seminal treatise \cite{Isaacs1965Differential}, Isaacs employed the principles of game theory, calculus of variations, and control theory, albeit unknown to him, to solve problems involving a dynamic conflict between multiple agents/players, \cite{Isaacs1965Differential}. Isaacs used the method of differential dynamic programming and introduced critical mathematical constructs such as dispersal, universal, and equivocal surfaces used to describe the optimal flow field in games and derive optimal saddle point strategies.

It is important to recognize some of the founders of static/dynamic games and optimal control including RAND scientists such as Richard E. Bellman, Leonard D. Berkovitz, David H. Blackwell, Melvin Dresher, Wendell H. Flemming, and John F. Nash. Early contributions on dynamic games in the former Soviet Union were published by N. Krasovskii, A. Melikyan, L. S. Pontryagin, and A. I. Subbotin \cite{Breitner2005Genesis}.  Furthermore, at the first conference dedicated exclusively to dynamic games: ``First International Conference on the Theory and Applications of Differential Games'' in Amherst, MA, September 29-October 1, 1969, organized by Ho and Leitmann. Isaacs, Berkovitz, Bernhardt, Blaquiere, Breakwell, Case, Friedman, Merz, Pontryagin, and Shubik were among the invited speakers. Although Pontraygin was not able to attend the meeting, these mathematicians and scientists planted the seeds of the theory of optimal control and differential games.

Of a large number of mathematicians and scientists which were involved in the development of differential games, Rufus Isaacs, Richard Bellman, John Breakwell, and Lev Pontryagin can be seen as principal contributors to the development of the theory of differential games; Isaacs being the father of Differential Games. His seminal work and his book highlighted the possible use of differential games to solve many problems, \cite{Isaacs1951Games, Isaacs1965Differential}. Bellman, known for the method of dynamic programming, provided a tool whereby state feedback optimal strategies could be directly obtained as opposed to methods based on necessary conditions as is the calculus of variations, \cite{Bellman1957, Lew1986}. Pontryagin, a Soviet mathematician, is recognized as developing his Maximum Principle (We'll refer to it as Pontryagin's Minimum Principle (PMP) for the rest of this paper) for necessary conditions for optimal control in the presence of hard constraints on control. Using methods derived by these four mathematicians, differential games were formulated and solved in many works to be described throughout this paper.

\section{Overview of Recent Results in Pursuit-Evasion}    \label{sec:litPE}
At the center of Differential Games lies the fundamental conflict of two parties known as ``Pursuit-Evasion''. Pursuit-evasion involves at least two agents or groups, labeled Pursuers and Evaders. The goal of a pursuer is to capture evading agents, while the converse is the goal of an evader, to avoid being captured by a pursuer. This is a zero-sum game where the cost/payoff is the time-to-capture. Basic questions arise, ``What path should an evader or pursuer take to achieve their goal of avoiding or ensuring capture; and, under optimal play, by either pursuer or evader, is capture at all possible?'' In this article, we begin by looking at this conflict and briefly discuss the current literature available describing strategies, methods, and applications.  We invite the interested reader to read the great historical and literature surveys of Isaacs' work documented in \cite{Breitner2005Genesis, Bernhard1998}.

The idea of pursuit-evasion differential games is not limited to physical entities chasing after one another; Isaacs defined kinematic equations which described the surfaces upon which states were constrained. Using these differential equations, one could propose problems in other research areas including but not limited to economics, sports, robotics, and air-combat. In this article, the focus is on differential games of pursuit-evasion, but it is important to note that the applications of these mathematical tools are not limited to these simple, toy problems, a syllogism for more complex scenarios which may not be suitable for the public domain.
\subsection{One Pursuer, One Evader (1v1)}    \label{sec:lit1v1}
The premise of a differential game starts with the conflict between two players who share a common performance functional; the goals of the pursuer and evader are counter to one another. These are minimax problems, that is, zero-sum-games, since an optimal solution for the performance functional is one in which the strategy of one of the players seeks to minimize the performance functional while the strategy of the other player tries to maximize it. Constraints on the player's come in the form of dynamics. These constraints can be linear or nonlinear. The classical problem of pursuit-evasion can be seen in an early work by Ho, Bryson, and Baron \cite{Ho1965Differential}. In their work, a two-player differential game was formulated in a Linear-Quadratic form to capture the basic pursuit-evasion conflict. Later, in the NASA technical report  \cite{Baron1966Differential}, differential models were employed in order to gauge the differences in performance between a manually piloted vehicle and an optimally controlled one as provided in the earlier work, \cite{Ho1965Differential}. The experiment showed that the use of differential games indeed provided useful information to pilots, but a cautionary statement at the end of the technical report stated that, "...differential game problems will, in general, be more complicated theoretically than their optimal control counterparts." The NASA report concluded that the idea of solving differential games was thought to be useful as information provided to a pilot, but not yet accepted to be a means of automatic control, a popular research topic today.

In a dissertation by Satimov \cite{Satimov1981a}, the application of differential games was envisioned for use in various fields such as economics and military operations. Satimov also stated that in the case of a single-player, differential games amount to optimal control problems, and that different modifications of Isaacs' method give the necessary and sufficient conditions for optimality. The relationship between optimal control and differential games is through the use of variational techniques, \cite{Isaacs1965Differential, Gelfand2000Calculus}. If all but one of the player's control laws are given, the differential game reverts to a one-sided optimal control problem.

\subsubsection{Homicidal Chauffeur Game}
In his seminal text, Isaacs proposed the famous ``Homicidal Chauffeur" problem \cite{Isaacs1965Differential}. In this game, a hypothetical slow but highly maneuverable holonomic pedestrian is pitted against a driver of a motor vehicle which is faster but less maneuverable (a.k.a. a Dubin's Car). In this somewhat macabre scenario, the driver attempts to run over the pedestrian. The question to be solved is: Under what circumstances, and with what strategy, can the driver of the car guarantee that he can always catch the pedestrian or conversely, the pedestrian guarantee that he can indefinitely elude the car, \cite{Isaacs1965Differential}. And, if the pedestrian's demise is guaranteed, what is the chauffeur's optimal strategy that will minimize the time-to-capture of the pedestrian, and what is the latter's strategy to maximize his time? This simple game is a classical 1-pursuer-1-evader problem used to represent military applications for reasons of acceptance and publication. Surveys have documented the history and notable work related to the ``Homicidal Chauffeur Differential Game,'' \cite{Marchal1975Analytical, Breitner2005Genesis, Patsko2009Homicidal, Falcone2005}, going into detail and expanding about the various aspects of this problem. 


A definitive work on the Homicidal Chauffeur Differential Game is Merz's Ph.D. thesis, \cite{Merz1971Homicidal}. Merz investigated the Homicidal Chauffeur Differential Game in great detail describing two new singular lines known as: ``switch envelope'' and ``focal line.'' These new lines further expand on Isaacs' ``barrier'', ``universal'', and ``equivocal'' singular lines. His work gives great detail and insight into the problem. Breakwell and Merz helped motivate the complete solution of the Homicidal Chauffeur game at a conference in 1969, \cite{Breakwell1969Complete}. Marchal also studied the Homicidal Chauffeur game in great detail describing how using Pontryagin's Minimum Principle could assist the interpretation of complex solutions, \cite{Marchal1975Analytical}.

\subsubsection{The Two Cars Differential Game}
\begin{figure}[htb]
\begin{center}
\includegraphics[width = 8.0cm]{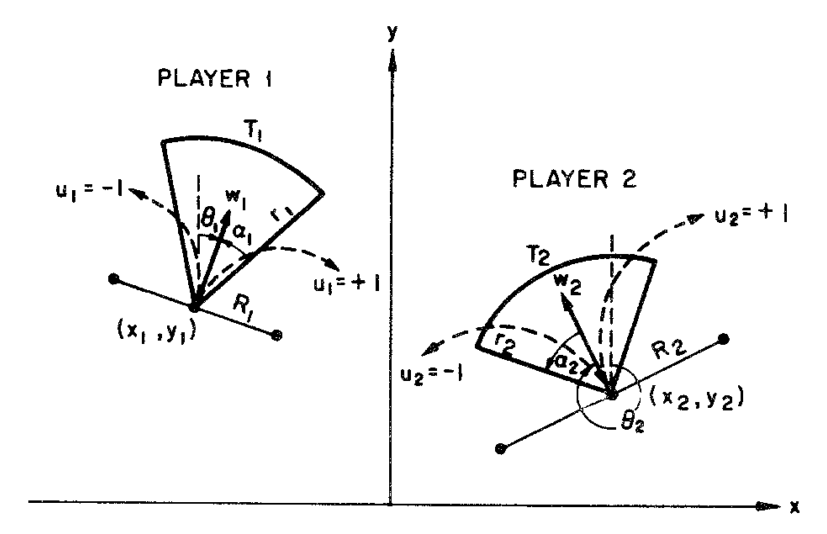}
\caption{The coordinate frames and turning radii used to describe the Two Cars Problem \cite{Getz1981Two} Radius $R_1$ and $R_2$ describe the minimum turning radii of each player and $u_1$ and $u_2$ describe the curves associated with a max rate turn. $w_1$ and $w_2$ describe each player's velocity.}
\label{fig:twoCars}
\end{center}
\end{figure}

A variation of the Homicidal Chauffeur Differential Game is the ``Two Cars'' Differential Game where two players, each controlling a car with minimum turning radius, are engaged in a pursuit-evasion game. In early work, Meier investigated the problem of Two Cars where both players had the same minimum turn radius, the pursuer was slower than the evader, and the capture was defined by coming inside the range, $l$, of the Evader, \cite{Meier1969New}. Another analysis of the Two-Car problem was performed by \cite{Getz1981Two, Getz1981Capturability}. In their papers, regions of capture, escape, and barrier surfaces between those regions were presented. Fig.~\ref{fig:twoCars} describes the geometry of the Two-Car Problem. Radius $R_1$ and $R_2$ describe the minimum turning radii of each player, $u_1$ and $u_2$ describe the curves associated with a max rate turn, and $w_1$ \& $w_2$ describe each player's velocity.

In \cite{Getz1981Two} both agents have sector based regions of capture, typical of an aerial dogfight; but, in \cite{Getz1981Capturability}, the regions of capture were different, describing a heterogeneous model of on-board weapon systems. Similarly, in \cite{Greenfeld1987a}, Greenfield looked into the two-car problem, endowing the pursuer with a surveillance capability of range, $l$. The objective, to escape the surveillance region in minimum time. In an earlier work, Lewin investigated a similar differential game called the ``Surveillance-Evasion'' Differential Game \cite{Lewin1975Surveillance}. In the game, the evader strives to escape as soon as possible from the pursuer's detection circle, while the pursuer's desire is the opposite. 

A complete analysis performed by Bera, Makkapati, and Kothari goes into detail of both games of kind and degree with studies on differing agent's speeds, capture radius and maneuverability constraints, \cite{Bera2017}. In their work, they develop the three-dimensional plots of the state space, highlighting the barrier and switching surfaces for the different scenarios.
\subsubsection{Pursuit-Evasion in Constrained Environments}
In an effort to consider differential games in a more realistic way, the introduction of boundaries and constraints on states allows for finite spaces and regions to be included in the game formulation. By imposing limitations on physical states, the pursuit-evasion differential game can be restricted to a bounded area or obstacles may be applied. A two player differential pursuit-evasion game where an obstacle is added to delay the pursuer or to avoid capture entirely was proposed by Fisac and Sastry \cite{Fisac2015a}. Similarly, the reference \cite{Oyler2016} considered pursuit-evasion games in the presence of obstacles which inhibit the motion of the players. In their work, the use of polygons, line-segments, and asymmetric obstacles (an obstacle which effects one player differently than another) are developed. Fuchs and Khargonekar motivated the use of escort regions through manipulation of the performance functional, \cite{Fuchs2017}. In \cite{Sundaram2017Pursuit} the pursuer and evader are constrained to road networks and, in\cite{Kalyanam2015}, a one-sided constrained is imposed where the pursuer is not restricted to a road network while the evader is. 
\subsubsection{Pursuit-Evasion with Incomplete Information}
In cases where one or more agents do not have complete information about the state of the game, these are problems called ``Differential Games of Incomplete Information.'' In his seminal work, \cite{Isaacs1965Differential} stated that the ability to pose problems which restricted information to the individual players ``...appears to be the most vital area for future research.'' Roxin and Tsokos, introduced a stochastic approach as a means of modeling partial information one agent might have relative to another, \cite{Roxin1969Definition}. Chernousko and Melikyan described a differential game where incomplete information is provided to one of the agents \cite{Chernousko1974Some}. This idea was proposed in order to account for information delay or gaps of information during game play. Yavin proposed an incomplete information pursuit-evasion differential game by restricting the pursuer's information on bearing and allowed the Evader to have perfect information in the engagement, \cite{Yavin1985}. Giovannangeli, Heymann and Rivlin tackled the problem of pursuit while avoiding convex obstacles by using Apollonius circles to provide paths in-which the pursuer's visibility of the evader is guaranteed throughout the engagement, \cite{Giovannangeli2010a}; the geometric concept of Apollonius circle will be formally defined in Section \ref{subsec:2p1Esol}. In \cite{Hexner1979Differential}, Hexner considered the problem where a parameter is unavailable to only one player at the beginning of the game, and the other has a probability density function describing that parameter was described. Pachter and Yavin investigated the effects of noise on the Homicidal Chauffeur problem by introducing stochasics to the pursuit-evasion differential game dynamics, \cite{Pachter1981Stochastic}. Battistini and Shima  also employed a stochastic model to overcome the limitations imposed by bearing-only measurements made by a pursuer, \cite{Battistini2013Bearings}. 

The authors of  \cite{Basimanebotlhe2014a} studied a differential game with nonlinear stochastic equation where two players are subjected to noisy measurements. The paper \cite{Lin2015Nash} presented an N-pursuer-1-evader differential where the evader can observe all the pursuers but the pursuers have limited observations of themselves and the evader. A pursuit-evasion game was studied in \cite{Kalyanam2016Pursuit} where a pursuer engages an evader using Unattended Ground Sensors (UGSs) that detect the evader's passage in a road network; when the pursuer arrives at an instrumented node, the UGSs inform the pursuer if and when the evader visited the node. 
\subsubsection{Pursuit Evasion in Aerial Engagements}
Applications of pursuit-evasion differential games relating to tactical air-to-air applications have been investigated. Shinar and Gutman developed a closed form solution to a 3-Dimensional missile-aircraft pursuit-evasion game, \cite{Shinar1978Recent}. Shinar also investigated a realistic pursuit-evasion engagement involving a missile engaged on an aircraft and air-to-air scenarios using variational methods, \cite{Shinar1980Realistic}. 
Considering naval applications, Pachter and Milch framed their two-player engagement as a Homicidal Chauffeur Differential Game where dynamics of the ships are taken into account,  \cite{Pachter1987Homicidal}. Greenwood designed a realistic differential game in 3-dimensions by modeling fighter aircraft, \cite{Greenwood1992a}. Greenwood used dynamics of two fighter aircraft in space and even considered firing envelopes as part of this barrier analysis. In \cite{Imado2005} a differential game is proposed which involves a pursuit-evasion engagement between a missile and an aircraft. In the game formulation, a nonlinear miss-distance was used as a payoff functional. In \cite{Shinar2007HybridPursuit} a pursuit-evasion game where the dynamics of the pursuer can be changed during the pursuit a finite number of times was investigated. In \cite{Shinar2009HybridEvade}, the evader has the ability to change their dynamics during the engagement a finite number of times. Merz investigated the problem of pursuit or evasion selection if both agents were endowed with capture sets and prior assignment had not been implemented, \cite{Merz1984}. Related to dog-fights and aerial combat, Merz's concern was with role assignment in pursuit-evasion differential games and of course the outcome. In more recent work \cite{Pachter2014Active,Garcia2015Active} the target area defense differential game has been investigated and it will be described in more detail in Section \ref{sec:litATDDG}.

\subsubsection{Other 1V1 Works}
%

Another example of posing the pursuit-evasion problems using simple motion kinematics in a differential game is in a work by Leitmann, \cite{Leitmann1968Simple}. In this paper, a simple differential game between a pursuer and evader was proposed, and variational techniques were applied to determine outcomes of the game where terminal miss distance was used as the payoff/cost functional.

In \cite{Calise1985Analysis}, Calise and Yu use simple motion kinematics as well as expanded control energy to formulate a game involving the pursuit-evasion of two aircraft at medium to long range. Using a reduced-order model with control energy, Calise and Yu are able to find trajectories similar to minimum time intercept using only four states to model the encounter.

The ``Lion and the Man'' differental game discussed by Quincampoix is a pursuit-evasion differential game where the lion pursues a man, \cite{Quincampoix2012Differential}. The lion and the man are free to change their velocity direction instantaneously but are limited to the intensity with which they do so. Since the lion is faster than the man, the regions of escape and capture are of interest and are numericaly determined.
\subsection{N Pursuers, 1 Evader (Nv1)} \label{sec:litnv1}
The two-pursuer-1-evader problem had been well documented, \cite{Hagedorn1976Differential, Pashkov1987a, Levchenkov1990a, Ganebny2011, Garcia2017a, Hayoun2017a}. Hagedorn and Breakwell investigated two pursuers engaging one evader which was required to pass between the two pursuers, \cite{Hagedorn1976Differential}. The game of degree was addressed in \cite{Pashkov1987a, Levchenkov1990a} by employing the following payoff functional: the distance between the object being pursued and the pursuer closest to it, when a fixed-time engagement terminates. The paper \cite{Ganebny2011} considered the three-player game in detail by briefly discussing the surfaces of the differential game when pursuers were both stronger, both weaker, and one stronger-one-weaker than the evader. The authors of \cite{Garcia2017a} studied the case where faster pursuers cooperate to capture a slower evader in minimum time. 
 Hayoun and Shima restrict the pursuer's controls to be bounded and their intercept times equal, \cite{Hayoun2017a}. Using two ``strong'' pursuers, closed-form optimal controls are derived, and it is shown that the addition of a second pursuer introduces a new singular zone to the game space in which the pursuers can guarantee equal misses, regardless of the evaders actions.

\begin{figure}[htb]
\begin{center}
\includegraphics[width=7.2cm]{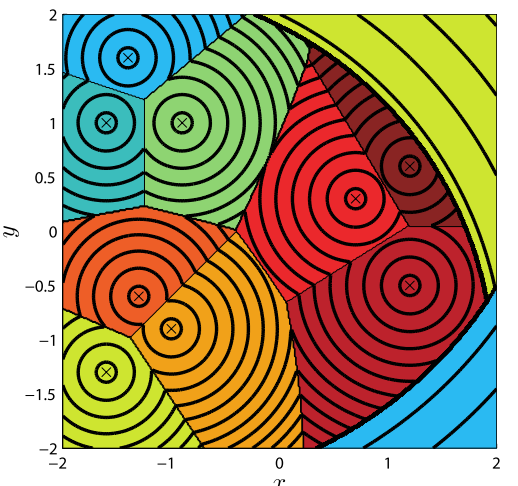}
\caption{A Voronoi Diagram \cite{Bakolas2010} describing the task allocation in an N-pursuer, 1-evader problem. Each color describes a cell from which a pursuer would capture an evader if starting in that cell. $\copyright$ 2010 IEEE}
	\label{fig:Voronoi}
\end{center}
\end{figure}

In the more general case, where there are N-pursuers and 1-evader, challenges with task allocation and strategy become more apparent. In order to aid the task allocation between the N-pursuers, Huang and Bakolas employ the Voronoi diagram  construct. It is used when capture of an evader within a bounded domain is considered, \cite{Huang2011,Bakolas2010}. In \cite{Borowko1985} sufficient conditions are presented for the existence of an evasion strategy where simple motion kinemtics for the players is considered. The paper \cite{Chodun1989} considered a more general case of N-pursuers engaged against one evader. Huang et al. employed a decentralized control scheme based on the Voronoi partition of the game domain, where the pursuers jointly shrink/minimize the area of the evader's Voronoi cell, \cite{Huang2011}. Fig.~\ref{fig:Voronoi} is a visualization of the individual pursuer cells from \cite{Bakolas2010}.

The N-pursuer-1-evader problem was investigated in \cite{Ibragimov2012} by focusing on guaranteed escape of the evader. 
The reference \cite{Kothari2017a} considered identical non-holonomic players where a computationally efficient algorithm to obtain approximate solutions is proposed. Solutions to a multi pursuer single-superior-evader pursuit-evasion differential game using fuzzy logic methods were developed in \cite{Evader2016,Al-Talabi2017}, where the formation control mechanism guarantees that the pursuers are distributed around the superior evader in order to avoid collision between pursuers and guarantees that the capture regions of each two adjacent pursuers are such that the capture of the fast evader is guaranteed. VonMoll et. al. also have considered multiple pursuers engaged on a single evader \cite{VonMoll2019MultiPursuer, VonMoll2018Evader, Pachter2019TwoV1}. By utilizing the Apollonius circles, they exploit the benefits of cooperation amongst the pursuers in order to reduce the capture time of the evader.

\subsection{1 Pursuer, M Evaders (1vM)} \label{sec:lit1vm}
The single pursuer against M-evaders differential game is a game where a pursuer tries to capture M-evaders in finite time. One challenge is to select the order in which the pursuer accomplishes his task in minimum time. In, these problems the pursuer is faster, more maneuverable, or has other advantages over the evaders.

The case which involves two evaders and one pursuer has received much attention, \cite{Fuchs2010a, Scott2013, Breakwell1979Point, Pachter1983One}. In \cite{Fuchs2010a} the case where a single pursuer engages two evaders was addressed. The goal of this work was to investigate a differential game where the pursuer tries to capture either of the evaders, minimizing its cost, and the evaders strive to escape the pursuer for as long as possible, increasing the payoff/functional of the pursuer. 
 Scott and Leonard investigated a scenario where two evaders employ coordinated strategies to evade a single pursuer, but also to keep them close to each other, \cite{Scott2013}. In \cite{Breakwell1979Point}, Breakwell and Hagedorn investigated the capture of two evaders in succession by one pursuer in minimum time. Pachter and Yavin  proposed a differential game of pursuit-evasion with one pursuer and two evaders, the motion of the players being affected by noise, \cite{Pachter1983One}. The stochastic game of degree is considered, where the pursuer strives to maximize the probability of his winning the game, i.e., of capturing at least one of the evaders. A 3-Dimensional pursuit-evasion differential game consisting of a pursuer engaged against a team consisting of two evaders was proposed in \cite{Abramyants1064Evasion}. The team of evaders consisted of a true evader and false decoy evader; the evaders coordinate their actions to ensure the true evader escapes without capture.

\begin{figure}[!htb]
\begin{center}
\includegraphics[width = 8.4cm]{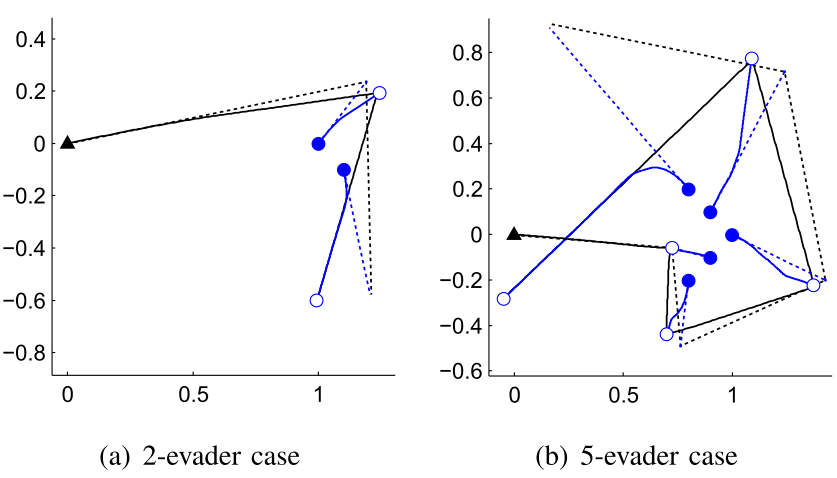}
\caption{Example of 1-pursuer-2-evader engagement and 1-pursuer-5-evader engagement found in \cite{liu2013evasion}. Since the pursuer is much more capable than the evaders capture of all agents is guaranteed. $\copyright$ 2013 IEEE}
\label{fig:1vm}
\end{center}
\end{figure}

In a more general case of one pursuer engaged against many evaders, the pursuer aims at capturing all evaders, while the evaders coordinate their escape. Liu and Zhou investigated a game involving a single-pursuer-multiple-evader pursuit-evasion game where a superior pursuer attempts to minimize the total capture time of all the evaders, \cite{liu2013evasion}. Fig.~\ref{fig:1vm} shows the capture of M-evaders in succession from \cite{liu2013evasion}. In \cite{Alias2015Simple} the study of capture time for many pursuers against one or more evaders is investigated. The authors of \cite{Wang2016a} formulate and solve a pursuit-evasion game in which a single faster player chases several homogenous evaders; a task allocation method to simulate the optimal engagements for ``fixed sequence capture'' and ``free sequence capture'' of the pursuer is applied in that reference. 
\subsection{N Pursuers, M Evaders (NvM)} \label{sec:litnvm}
The most general case of N-pursuers and M-evaders allows for more complex engagements to be analyzed. Katz et. al investigated a zero-sum differential game formulation for the control of military air operations using the method of characteristics, \cite{Katz2005a}. Although their examples were shown for 1-Pursuer 1-Evader, their work has extensions to N-pursuer n-Evader problems. Rusnak proposed a dynamic game called ``The Lady and the Body-Guards versus the Bandits'', \cite{Rusnak2005Lady}. The Bandits' team objective is to capture the Lady while the Lady and her Body-Guards objective is to prevent it. The Body-Guards are trying to intercept the Bandits prior to their arrival to the proximity of the Lady. The formulation and solution of the game is presented. As described in Section \ref{sec:litATDDG}, this problem is a similar problem, but with more players involved. A creative approach to handing the task allocation of many agents was proposed by Bakolas and Tsiotras by employing the Voronoi diagram construct, \cite{Bakolas2010}. Using the Voronoi diagram, such that a pursuer residing inside a given set of the partition can intercept a moving target faster than any other pursuer outside the set. Another means of task allocation was proposed in \cite{Awheda2016a} where Awheda and Schwartz proposed a fuzzy logic based decentralized control scheme using the Apollonius circles construct. 
A formal analysis of both optimal guidance and optimal assignments of $N$ pursuers to $M$ evaders, for $N\geq M$, was undertaken in \cite{garcia2019multiple}. 
Finally, the papers  \cite{margellos2011hamilton,Chen2017MultiplayerOutcomes,zhou2018efficient} address multiplayer differential games based on numerical solutions of Hamilton-Jacobi-Isaacs (HJI) equations.

\section{Two Cutters and Fugitive Ship Differential Game} \label{sec:TwoCutters}
This section provides a detailed treatment of Isaacs' classical problem of two cutters and fugitive ship differential game where two faster pursuers cooperate to capture a slower evader in minimum time. The evader, knowing that is being pursued by two cooperative and fast pursuers, tries to maximize the capture time. Isaacs proposed the players' optimal strategies in \cite{Isaacs1965Differential} and these strategies were verified in \cite{Garcia2017a} where the Value function was obtained and it was shown to be continuous, continuously differentiable, or $\mathcal{C}^1$, and the solution of the Hamilton-Jacobi-Isaacs (HJI) Partial Differential Equation (PDE). This section extends the presentation in \cite{Garcia2017a} by offering a complete proof of the verification theorem pertaining to the two cutters and fugitive ship differential game, it also identifies the game's singular surface which is of dispersal type, and it illustrates the applicability of the solution of this differential game to address interesting and more complex pursuit-evasion games with multiple pursuers and multiple evaders.

One of the fundamental concepts in differential games is the HJI PDE which provides a sufficient condition for the existence of a Nash equilibrium. Regardless of how the Value function is obtained, the HJI equation can be used to verify that such function indeed meets the optimality conditions.
Moreover, since the gradient of the Value function takes the place of the co-states in the PMP, a direct solution of the HJI PDE can be used to synthesize the optimal strategies. In this section we use the former objective, that is, verification. In the following, we summarize the sufficient conditions for a Nash equilibrium.

Consider the following dynamics
\begin{align}
 \left.
	\begin{array}{l l}
	\dot{\textbf{x}}=f(t,\textbf{x}(t),u(t),v(t))
	\end{array}  \right.   \label{eq:VerSumDyn}
\end{align}
for $t\in[0,t_f]$ and $\textbf{x}(0)=\textbf{x}_0$. The input variables $u(t)$ and $v(t)$ are the controls of each team, which are required to belong to an appropriate subspace. The performance functional is
\begin{align}
 \left.
	\begin{array}{l l}
	J=q(\textbf{x}(t_f)) +\int_0^{t_f} g(t,\textbf{x}(t),u(t),v(t)) dt
	\end{array}  \right.   \label{eq:VerSumCost}
\end{align}
State feedback policies are considered, that is, $u(t)=\mu(t,\textbf{x}(t))$ and $v(t)=\nu(t,\textbf{x}(t))$, for $t\in[0,t_f]$. 
\begin{theorem}
\cite{hespanha2017}. Assume that there exists a $\mathcal{C}^1$ function $V(t,\textbf{x}(t))$ that satisfies the HJI PDE  
\begin{align}
 \left.
	\begin{array}{l l}
-\frac{\partial V}{\partial t} =\frac{\partial V}{\partial \textbf{x}}\cdot  \textbf{f}(t,\textbf{x},u^*,v^*) + g(t,\textbf{x},u^*,v^*)
\end{array}  \right.   \label{eq:HJIPDE}
\end{align}
with $V(t_f,\textbf{x})=q(\textbf{x})$. Then, the pair of strategies $u^*=\mu^*(t,\textbf{x}(t))$ and $v^*=\nu^*(t,\textbf{x}(t))$ is a saddle-point equilibrium in state feedback policies.
\end{theorem}
\ \\
This result will be used to verify the solution of the two cutters and fugitive ship differential game in Section \ref{subsec:Verification} but first we provide a formulation of the problem.

\subsection{Problem Formulation} \label{sec:2PEproblem}

Consider one Evader $E$ and two Pursuers $P_1$ and $P_2$ with simple motion dynamics and constant speeds $v_E$, $v_{P1}$, and $v_{P2}$, respectively; the players are holonomic. It is assumed that the pursuers are faster than the evader.
Through explicit cooperation, the Pursuers aim at intercepting the Evader in minimum time, while the Evader tries to maximize the capture time. The advantage of having two cooperative Pursuers is that, depending on the initial positions of the players, the cooperative capture time is less than the capture time without cooperation.
This problem has important repercussions; although one pursuer is guaranteed to capture a slower Evader in an open plane, there exist situations where the Evader may win the game by reaching a certain subset of the space or exiting the domain of the game. For instance, consider an evading ship fleeing a pursuer and trying to exit the influence zone of the later. In this case the Pursuer would benefit from a second Pursuer and cooperation that strive to minimize the capture time preventing the Evader from reaching a safe haven. 

The controls of $E$, $P_1$, and $P_2$ are their respective instantaneous headings $\phi$, $\psi_1$, and $\psi_2$. The states of $E$, $P_1$, and $P_2$ are defined by their Cartesian coordinates in the realistic plane $\textbf{x}_E=(x_E,y_E)$, $\textbf{x}_{P1}=(x_{P1},y_{P1})$, and $\textbf{x}_{P2}=(x_{P2},y_{P2})$, respectively. The complete state of the game is defined by $\textbf{x}:= (x_{E}, y_{E}, x_{P1}, y_{P1}, x_{P2}, y_{P2})\in \mathbb{R}^6$. The Evader's control variable is his instantaneous heading angle, $\textbf{u}_E=\left\{\phi\right\}$. $P_1$ and $P_2$ affect the state of the game by choosing the instantaneous respective headings, $\psi_1$ and $\psi_2$, so the Pursuers' control variable is $\textbf{u}_P=\left\{\psi_1,\psi_2\right\}$. The normalized dynamics in the realistic plane $\dot{\textbf{x}}=\textbf{f}(\textbf{x},\textbf{u}_P,\textbf{u}_E)$ are specified by the system of ordinary differential equations
\begin{align}
 \left.
	\begin{array}{l l}
	\dot{x}_E&=\cos\phi,  \ \ \  \ \ \ \ \   
	\dot{y}_E=\sin\phi       \\
  \dot{x}_{P1}&=\beta_1\cos\psi_1,   \ \ \    
	\dot{y}_{P1}=\beta_1\sin\psi_1   \\
	\dot{x}_{P2}&=\beta_2\cos\psi_2, \ \ \    
	\dot{y}_{P2}=\beta_2\sin\psi_2   
	\end{array}  \right.   \label{eq:xT}
\end{align}
where $\beta_1=\frac{v_{P1}}{v_E}>1$ and $\beta_2=\frac{v_{P2}}{v_E}>1$ are the normalized speeds of the two Pursuers and they are not necessarily equal. The initial state is $\textbf{x}_0 := (x_{E_0}, y_{E_0}, x_{P1_0}, y_{P1_0}, x_{P2_0}, y_{P2_0}) = \textbf{x}(t_0)$.
We confine our attention to point capture, so the game terminates at time $t_f$ when the state of the system satisfies the terminal condition
\begin{align}
 \left.
	\begin{array}{l l}
  &a) \ \ x_{P1}=x_E, \ \ \  	 y_{P1}=y_E,	  or   \\
	&b) \ \ x_{P2}=x_E, \ \ \  	 y_{P2}=y_E,	   or  \\
	&c) \ \ both \ a) \ and \ b) \ apply.  
\end{array}  \right.  \label{eq:PC}
\end{align}
The terminal time $t_f$ is defined as the time instant when the state of the system satisfies any one of the conditions in \eqref{eq:PC}, the terminal state being $\textbf{x}_f := (x_{E_f}, y_{E_f}, x_{P1_f}, y_{P1_f}, x_{P2_f}, y_{P2_f}) = \textbf{x}(t_f)$. 
The Evader strives to maximize the capture time while the Pursuers cooperate to minimize the capture time, so the performance functional is
	\begin{align}
	\min_{\psi_1,\psi_2} \max_\phi  \int_0^{t_f}dt   \label{eq:mmJ}
\end{align}
subject to \eqref{eq:xT}-\eqref{eq:PC}.

\subsection{Solution}   \label{subsec:2p1Esol}
In this section we follow the procedure to derive regular solutions of differential games \cite{Basar1998Dynamic}, \cite{Lewin94}. 
We introduce the co-states $ \lambda:=(\lambda_{x_E},\lambda_{y_E}, \lambda_{x_{P1}}, \lambda_{y_{P1}}, \lambda_{x_{P2}}, \lambda_{y_{P2}}) \in \mathbb{R}^6$ and the Hamiltonian of the differential game is formulated as follows
\begin{align}
  \left.
	\begin{array}{l l}
	\mathcal{H}&=1+\lambda_{x_E}\cos\phi + \lambda_{y_E}\sin\phi  \\
	&~~+ \beta_1\lambda_{x_{P1}}\cos\psi_1 + \beta_1\lambda_{y_{P1}}\sin\psi_1 \\
	&~~+ \beta_2\lambda_{x_{P2}}\cos\psi_2 + \beta_2\lambda_{y_{P2}}\sin\psi_2. 
\end{array}  \right.   \label{eq:Ham}
\end{align}
The optimal control inputs (in terms of the co-state variables) are obtained from 
\begin{align}
    \min_{\psi_i,\psi_2} \max_\phi \mathcal{H} \equiv 0 	\label{eq:minmax}
\end{align}
and they are given by
\begin{align}
\left.
	\begin{array}{l l}
   \cos\phi^*&={\scriptstyle \frac{\lambda_{x_E}}{\sqrt{\lambda_{x_E}^2+\lambda_{y_E}^2}}}, \ \ \ \ \ 
	\sin\phi^*={\scriptstyle \frac{\lambda_{y_E}}{\sqrt{\lambda_{x_E}^2+\lambda_{y_E}^2}}}
	\end{array}  \right. \label{eq:chico}  
\end{align}
\begin{align}
  \left.
	\begin{array}{l l}
	\cos\psi_1^*={\scriptstyle -\frac{\lambda_{x_{P1}}}{\sqrt{\lambda_{x_{P1}}^2+\lambda_{y_{P1}}^2}}}, \ \ \
	\sin\psi_1^*={\scriptstyle -\frac{\lambda_{y_{P1}}}{\sqrt{\lambda_{x_{P1}}^2+\lambda_{y_{P1}}^2}}} 
	\end{array}  \right. \label{eq:psico} 
\end{align}
\begin{align}
	\left.
	\begin{array}{l l}
	  \cos\psi_2^*={\scriptstyle -\frac{\lambda_{x_{P2}}}{\sqrt{\lambda_{x_{P2}}^2+\lambda_{y_{P2}}^2}}},  \ \ \ 
		\sin\psi_2^*={\scriptstyle -\frac{\lambda_{y_{P2}}}{\sqrt{\lambda_{x_{P2}}^2+\lambda_{y_{P2}}^2}}}. 
	\end{array}  \right. \label{eq:phico} 
\end{align}
Additionally, we can determine the co-state dynamics by calculating $\dot{\lambda}=-\frac{\partial H}{\partial \textbf{x}}$. Hence, the co-state dynamics are: $\dot{\lambda}_{x_E}=\dot{\lambda}_{y_E}=\dot{\lambda}_{x_{P1}}=\dot{\lambda}_{y_{P1}}=\dot{\lambda}_{x_{P2}}=\dot{\lambda}_{y_{P1}}=0$; then, all co-states are constant and the optimal control inputs \eqref{eq:chico}-\eqref{eq:phico} are constant as well. 
Because the optimal headings are constant then the players' optimal trajectories are straight lines. Also, since the terminal condition is point capture \eqref{eq:PC} we have the following: For case a) in \eqref{eq:PC}:
\begin{align}
\left.
	\begin{array}{l l}
   \!\! x_{P1_f}\! =\! x_{P1_0}\!+\!\beta_1t_{f_1}\!\cos\psi_1^* \!=\! x_{E_0} \!+\! t_{f_1}\cos\phi^* \!=\! x_{E_f}  \\
	\!\! y_{P1_f} \!=\! y_{P1_0}\!+\!\beta_1t_{f_1}\sin\psi_1^* = y_{E_0} \!+\! t_{f_1}\sin\phi^* \!=\! y_{E_f}
\end{array}  \right.   \label{eq:P1eqs}
\end{align}
where $t_{f_1}\geq 0$ can be written in terms of $\phi^*$ as follows
\begin{align}
  t_{f_1}=c_1\cos(\phi^*-\lambda_1) + \sqrt{c_1^2\cos^2(\phi^*-\lambda_1) + c_1r_{1_0}}  \label{eq:tf1}
\end{align}
where 
\begin{align}
\lambda_1&=\arctan\Big(\frac{y_{E_0}-y_{P1_0}}{x_{E_0}-x_{P1_0}}\Big) \nonumber \\
r_{1_0}&=\sqrt{(x_{E_0}-x_{P1_0})^2+(y_{E_0}-y_{P1_0})^2}  \nonumber  \\
c_1&=\frac{1}{\beta_1^2-1}r_{1_0}.  \nonumber
\end{align}
For case b) in \eqref{eq:PC} we have:
\begin{align}
\left.
	\begin{array}{l l}
  x_{P2_0}+\beta_2t_{f_2}\cos\psi_2^* = x_{E_0} + t_{f_2}\cos\phi^*   \\
	y_{P2_0}+\beta_2t_{f_2}\sin\psi_2^* = y_{E_0} + t_{f_2}\sin\phi^* 
\end{array}  \right.   \label{eq:P2eqs}
\end{align}
where $t_{f_2}\geq 0$ can be written in terms of $\phi^*$ as follows
\begin{align}
  t_{f_2}=c_2\cos(\phi^*-\lambda_2) + \sqrt{c_2^2\cos^2(\phi^*-\lambda_2) + c_2r_{2_0}}   \label{eq:tf2}
\end{align}
where 
\begin{align}
\lambda_2&=\arctan\Big(\frac{y_{E_0}-y_{P2_0}}{x_{E_0}-x_{P2_0}}\Big) \nonumber \\
r_{2_0}&=\sqrt{(x_{E_0}-x_{P2_0})^2+(y_{E_0}-y_{P2_0})^2}  \nonumber \\
c_2&=\frac{1}{\beta_2^2-1}r_{2_0}.  \nonumber
\end{align}
Finally, for the most interesting case c) in \eqref{eq:PC}, both \eqref{eq:P1eqs} and \eqref{eq:P2eqs} are simultaneously satisfied for the same heading $\phi^*$.

Because $\beta_1>1$ (respectively $\beta_2>1$) there exists a (unique) solution $\psi_1$ (respectively $\psi_2$) to equation \eqref{eq:P1eqs} (respectively \eqref{eq:P2eqs}) for any $\phi\in[0,2\pi]$. However, the opposite is not true. There exists a range of heading values for $\psi_1$ (respectively $\psi_2$) such that no solution $\phi$ exists (these Pursuer headings correspond to the cases where the Pursuer runs away from the Evader). Note that $t_{f_1}=0$ (respectively $t_{f_2}=0$) only if $r_1=0$ (respectively, $r_2=0$), where $r_i=\sqrt{(x_E-x_{Pi})^2+(y_E-y_{Pi})^2}$, for $i=1,2$. In such a case the Evader has been captured and the Game has ended.

Let us assume that condition a) in \eqref{eq:PC} is active, whereupon the optimal Evader strategy is obtained by solving
\begin{align}
\left.
	\begin{array}{l l}
  \max_\phi t_{f_1}   \\
	\text{subject to eq. \eqref{eq:P1eqs}}
\end{array}  \right.   
\end{align}
which can be computed by solving for $\phi$ in
\begin{align}
\left.
	\begin{array}{l l}
  &\frac{dt_{f_1}}{d\phi}=-c_1\sin(\phi-\lambda_1)  \\
	&~~~~~~-\frac{c_1^2\cos(\phi^*-\lambda_1)\sin(\phi^*-\lambda_1)}{\sqrt{c_1^2\cos^2(\phi^*-\lambda_1) + c_1r_{1_0}}} = 0  \\
	\Rightarrow   & \sin(\phi-\lambda_1) = 0  \ \
	\Rightarrow \ \phi^*=\lambda_1    
\end{array}  \right.   \label{eq:Opt1}
\end{align}
as expected (the heading $\phi=\lambda_1+\pi$ will minimize $t_{f_1}$ and since $E$ chooses its heading $\phi$ then it would choose the one that maximizes $t_{f_1}$).

Similarly, if condition b) in \eqref{eq:PC} is active, then the optimal Evader strategy is $\phi^*=\lambda_2$, as expected. In cases a) or b) the solution of the Two-Pursuer One-Evader Differential Game (2P1EDG) simplifies to the solution of the differential game with one pursuer and one Evader where the Evader runs directly away from the pursuer. 

Since it is not evident at the beginning which condition, a), b), or c), is active, then the individual conditions a) and b) are checked first as follows:
Consider $\phi=\lambda_1$ and obtain $t_{f_1}(\lambda_1)$ using \eqref{eq:tf1}. Then, substitute $\phi=\lambda_1$ in \eqref{eq:tf2} to compute $t_{f_2}(\lambda_1)$ and compare: 
\begin{align}
  \text{if} \ t_{f_1}(\lambda_1) \leq t_{f_2}(\lambda_1)   \label{eq:ca}
\end{align}
then $\phi^*=\lambda_1$ and there is nothing Pursuer 2 can do to decrease the capture time $t_{f_1}(\phi^*=\lambda_1)$.

Also, consider $\phi=\lambda_2$ and obtain $t_{f_2}(\lambda_2)$ using \eqref{eq:tf2}. Then, substitute $\phi=\lambda_2$ in \eqref{eq:tf1} to compute $t_{f_1}(\lambda_2)$ and compare: 
\begin{align}
  \text{if} \ t_{f_2}(\lambda_2) \leq t_{f_1}(\lambda_2)   \label{eq:cb}
\end{align}
then $\phi^*=\lambda_2$ and there is nothing Pursuer 1 can do to decrease the capture time $t_{f_2}(\phi^*=\lambda_2)$.

If neither, condition \eqref{eq:ca} nor condition \eqref{eq:cb}, is satisfied, we check if the interesting condition c) is active. In such a case of simultaneous capture of $E$ by $P_1$ and $P_2$, both equations \eqref{eq:P1eqs} and \eqref{eq:P2eqs} need to hold for the same $\phi^*$ and for the same $t_f=t_{f_1}=t_{f_2}$, where $t_{f_1}$ and $t_{f_2}$ are given by \eqref{eq:tf1} and \eqref{eq:tf2}, respectively. Additionally, for initial conditions such that c) is active, the capture time will be reduced compared to the individual cases where the game only involves $E$ and $P_1$ or only $E$ and $P_2$. Hence, the complete solution of the 2P1EDG should define the optimal saddle point strategies of every player and the outcome of the game in terms of which pursuer will actually capture the Evader. 

In case c) in \eqref{eq:PC} the triple $(\phi,\psi_1,\psi_2)$ which simultaneously satisfies \eqref{eq:P1eqs}-\eqref{eq:tf2} with $t_{f_1}=t_{f_2}$ can be obtained from the intersection of two Apollonius circles. The first circle is based on $r_{1_0}$ and the speed ratio $\beta_1$. The second circle is based on $r_{2_0}$ and the speed ratio $\beta_2$. 
The center of the Apollonius circle between $E$ and $P_i$, for $i=1,2$, is at a distance $c_i=\frac{1}{\beta_i^2-1}r_{i_0}$ from $E$ and its radius is $\rho_i=\frac{\beta_i}{\beta_i^2-1}r_{i_0}$. The coordinates of the circle center are 
\begin{align}
   a_i=x_E + \rho_i \cos\lambda_i, \ \ \ 
		 b_i=y_E + \rho_i \sin\lambda_i	   \nonumber
\end{align}
and the equation of the circle is 
\begin{align}
(x-a_i)^2+(y-b_i)^2=\rho_i^2  \label{eq:ApCircle}
\end{align}
for $i=1,2$. From the two points of intersection of the circles, the Evader chooses the point that maximizes the capture time, as expected; this point is the intersection point of the Apollonius circles \eqref{eq:ApCircle} which is the farthest away from the Evader's position $(x_E,y_E)$. Let $I:(x_I,y_I)$ be the aimpoint and we have that $\phi_s^*=\arctan(\frac{y_I-y_E}{x_I-x_E})$, $\psi_{1s}^*=\arctan(\frac{y_I-y_{P1}}{x_I-x_{P_1}})$, and $\psi_{2s}^*=\arctan(\frac{y_I-y_{P2}}{x_I-x_{P2}})$ which guarantee simultaneous interception by both Pursuers. In the case of simultaneous capture we have that $\lambda_m \leq \phi_s^* \leq \lambda_M$, where $\lambda_m=\min\left\{\lambda_1,\lambda_2\right\}$ and $\lambda_M=\max\left\{\lambda_1,\lambda_2\right\}$.

In this section we described a method that provides a simple way to design a state feedback pursuit strategy for the 2P1EDG and it is summarized as follows.

\begin{proposition} \label{prop:method}
Consider the Two-Pursuer One-Evader Differential Game \eqref{eq:xT}-\eqref{eq:mmJ} where $\beta_1>1$ and $\beta_2>1$. The solution of the differential game is: 
\begin{align}
	\left.
	\begin{array}{l l}
	     \phi^*=\lambda_1, \ \psi_1^*=\lambda_1 \ \ \  \ \ \ \  \ \ \ \ \ \ \ \ \ \  \ \text{if} \ t_{f_1}(\lambda_1)\leq t_{f_2}(\lambda_1)  \\
				\phi^*=\lambda_2, \ \psi_2^*=\lambda_2 \ \ \ \ \ \ \ \ \	\ \ \ \ \ \ \ \ \ \text{if}  \ t_{f_2}(\lambda_2)\leq t_{f_1}(\lambda_2)  \\
	     \phi^*=\phi_s^*, \ \psi_1^*= \psi_{1s}^*, \ \psi_2^*= \psi_{2s}^*  \ \ \   \text{otherwise}
	\end{array}  \right.  \label{eq:DGsol}
\end{align}
\end{proposition}


\subsection{Verification}   \label{subsec:Verification}
In this section we obtain the Value function and it is verified that the Value function is continuous and continuously differentiable and that it also satisfies the Hamilton-Jacobi-Isaacs (HJI) Partial Differential Equation (PDE). 

Based on the different outcomes of the Differential Game, the capture set, namely, the whole state space, is partitioned into three subsets. These three subsets are: 

- $\mathcal{R}_1$: where $E$ is only captured by $P_1$.

- $\mathcal{R}_2$: where $E$ is only captured by $P_2$.

- $\mathcal{R}_s$: where $E$ is simultaneously captured by $P_1$ and $P_2$.

\begin{figure}
	\begin{center}
		\includegraphics[width=8.4cm,height=6.5cm,trim=.9cm .5cm .8cm .2cm]{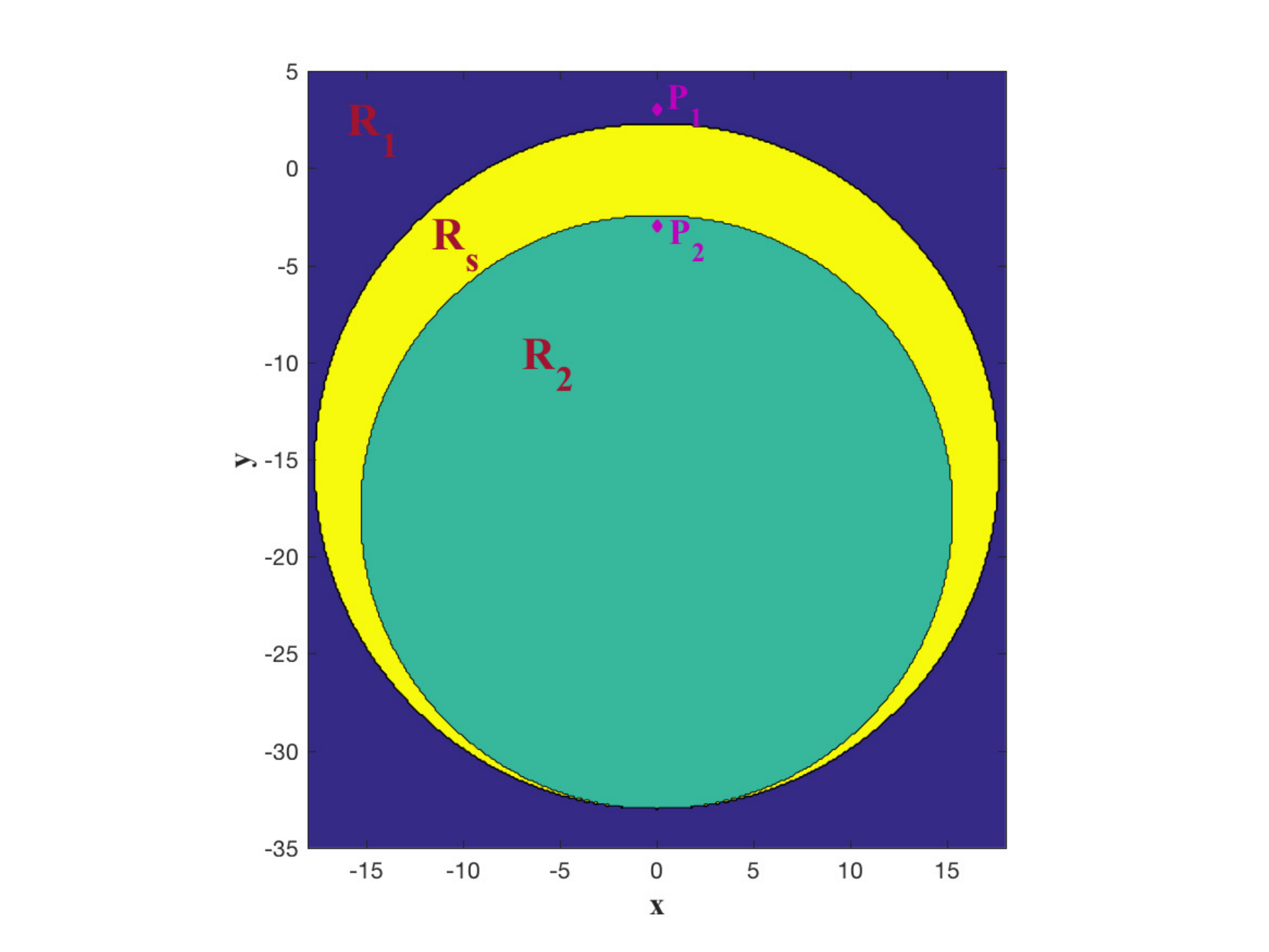}
	\caption{Regions $\mathcal{R}_1$, $\mathcal{R}_2$, and $\mathcal{R}_s$. The initial position of the Pursuers are $P_1:(0,3)$ and $P_2:(0,-3)$ and the speed ratio parameters are $\beta_1=1.3$ and $\beta_2=1.25$. 
}
	\label{fig:Regs}
	\end{center}
\end{figure}

In addition, for given speed ratio parameters $\beta_1$ and $\beta_2$ the subsets can be characterized using the conditions in \eqref{eq:ca} and \eqref{eq:cb}. 
The subset of states that result in \eqref{eq:ca} (respectively, \eqref{eq:cb}) holding with equality represent the boundary between regions $\mathcal{R}_1$ and $\mathcal{R}_s$ (respectively, regions $\mathcal{R}_2$ and $\mathcal{R}_s$). An example of these regions is shown in Fig. \ref{fig:Regs}.
Provided that every agent plays optimally, Fig. \ref{fig:Regs} represents the property that if $E$ is initially located inside $\mathcal{R}_1$, then the outcome of the game is that $E$ is captured only by $P_1$. If $E$ is in $\mathcal{R}_2$, it is captured only by $P_2$ and if $E$ is in $\mathcal{R}_s$, then it is captured simultaneously by both Pursuers. 

We wish to obtain a state feedback solution. At any time $t$, $0<t<t_f$ the state is $\textbf{x}(t)=[x_E(t),y_E(t),x_{P_1}(t),$ $y_{P_1}(t),x_{P_2}(t),y_{P_2}(t)]^T$.
The  Value function is given by 
\begin{align}
	V(\textbf{x})&=\left\{
	\begin{array}{l l}
	\frac{\sqrt{(x_E-x_{P1})^2+(y_E-y_{P1})^2}}{\beta_1-1} \qquad   \forall \ \textbf{x}\in \mathcal{R}_1   \\    
	\frac{\sqrt{(x_E-x_{P2})^2+(y_E-y_{P2})^2}}{\beta_2-1}\qquad   \forall \ \textbf{x}\in \mathcal{R}_2    \\    
	\frac{F_1(\textbf{x})t_{f_2}(\textbf{x})-F_2(\textbf{x})t_{f_1}(\textbf{x})}{F_1(\textbf{x})-F_2(\textbf{x})}  \qquad \ \ \  	\forall \ \textbf{x}\in \mathcal{R}_s  
	\end{array}  \right.  \label{eq:ValueFn}
\end{align}
 where 
\begin{align}
\left.
	 \begin{array}{l l} 
t_{f_i}(\textbf{x})\!\!=\!{\scriptstyle \frac{1}{\beta_i^2-1}} \big[(x_E\!-\!x_{P_i})\cos\phi^* \!+\! (y_E\!-\!y_{P_i})\sin\phi^*  \\
	 \qquad \ \  + \big([(x_E\!-\!x_{P_i})\cos\phi^* + (y_E\!-\!y_{P_i})\sin\phi^*]^2  \\
	\qquad \ \ \  +(\beta_i^2\!-\!1)[(x_E\!-\!x_{P_i})^2 \!+\! (y_E\!-\!y_{P_i})^2] \big)^{1/2}  \big]
\end{array}   \right.   \label{eq:tfi}
\end{align}
\begin{align}
\left.
	 \begin{array}{l l} 
F_{i}(\textbf{x})\!\!=\big[(y_E-y_{P_i})\cos\phi^* - (x_E-x_{P_i})\sin\phi^* \big]  \\
	 \qquad \ \ \times \big([(x_E\!-\!x_{P_i})\cos\phi^* + (y_E-y_{P_i})\sin\phi^*]^2  \\
	  \qquad \ \ \ +(\beta_i^2\!-\!1)[(x_E\!-\!x_{P_i})^2 \!+\! (y_E\!-\!y_{P_i})^2] \big)^{-1/2}    
\end{array}   \right.   \label{eq:Fi}
\end{align}
and $\beta_i>1$, for $i=1,2$, where $\phi^*(\textbf{x},\beta_1,\beta_2)$ is obtained from the points of intersection of the two Apollonius circles. 
The optimal Evader heading $\phi^*$ is \textit{only a function} of the state $\textbf{x}$ (the positions of the three players) and of the speed ratio parameters $\beta_1$ and $\beta_2$. Due to the complexity of an explicit expression of $\phi^*$, the terms $\cos\phi^*$ and $\sin\phi^*$ remain in $V(\textbf{x})$. The physical meaning of the Value function is the capture time under optimal play. The Value function takes different forms depending on which region the Evader is located; in any case, it should be a function of the state of the system. In region $ \mathcal{R}_1$, it is only a function of the states of $E$ and $P_1$. In region $ \mathcal{R}_2$, it is only a function of the states of $E$ and $P_2$. Finally, in region $ \mathcal{R}_s$, $V$ has to be correctly written in terms of the states of $E$, $P_1$, and $P_2$.

The last form of $V(\textbf{x})$ in \eqref{eq:ValueFn}, where $\textbf{x}\in \mathcal{R}_s$, is a convex combination of the Value of the game, the capture time, using the states of both Pursuers. 
When $\textbf{x}\in \mathcal{R}_s$ the condition for simultaneous capture is that $t_{f_1}(\phi^*)=t_{f_2}(\phi^*)$. In such a case, we define $t_f(\phi^*)\equiv t_{f_1}(\phi^*)=t_{f_2}(\phi^*)$.

\textit{Dispersal surface}. A dispersal surface $\mathcal{D}\in \mathcal{R}_s$ exists if the two intersections of the Apollonius circles are located at the same distance from the Evader's position; hence, two optimal strategies exist. This dispersal surface will be illustrated in Section \ref{subsec:Dispersal}. Let us define $\mathcal{R}_s^-=\mathcal{R}_s - \mathcal{D}$ and define the subset of the state space where regular solutions apply $\mathcal{R}=\mathcal{R}_1\cup\mathcal{R}_2\cup\mathcal{R}_s^-$.

\begin{theorem}
Consider the Two-Pursuer One-Evader Differential Game \eqref{eq:xT}-\eqref{eq:mmJ}. The solution of the differential game is given by Proposition \ref{prop:method}. The corresponding Value function is given by \eqref{eq:ValueFn}. Then, the Value function \eqref{eq:ValueFn} is continuous and continuously differentiable over $\mathcal{R}$ and it satisfies the HJI equation for any $\textbf{x}\in\mathcal{R}$.
\end{theorem}
\textit{Proof}. We start by analyzing the gradient of $V(\textbf{x})$ in regions $\mathcal{R}_1$ and $\mathcal{R}_2$. In the case where $\textbf{x}\in \mathcal{R}_1$ we can write 
\begin{align}
\left.
	 \begin{array}{l l} 
	 \frac{\partial V}{\partial x_E} & = \frac{1}{\beta_1-1}\cos\lambda_1, \ \  	
	 \frac{\partial V}{\partial y_E} = 	\frac{1}{\beta_1-1}\sin\lambda_1 \\
	 \frac{\partial V}{\partial x_{P1}} &= -\frac{\partial V}{\partial x_E}, \ \ \ \ \ \ \ \ \
 	 \frac{\partial V}{\partial y_{P1}} = -\frac{\partial V}{\partial y_E}	
\end{array}   \right.   \label{eq:Grad1}
\end{align}
where 
 \begin{align}
\left.
	 \begin{array}{l l} 
	\cos\lambda_1&=\frac{x_E-x_{P1}}{\sqrt{(x_E-x_{P1})^2+(y_E-y_{P1})^2}}  \\
	\sin\lambda_1&=\frac{y_E-y_{P1}}{\sqrt{(x_E-x_{P1})^2+(y_E-y_{P1})^2}} .
\end{array}   \right.  \nonumber
\end{align}
Therefore, the Value function $V(\textbf{x})$ is $C^1$ inside region $\mathcal{R}_1$. Similar expressions to \eqref{eq:Grad1} can be obtained when $\textbf{x}\in \mathcal{R}_2$ to show that $V(\textbf{x})$ is $C^1$ inside region $\mathcal{R}_2$. 

Now, the partial derivative of the Value function with respect to each element of the state $\textbf{x}$ when $\textbf{x}\in \mathcal{R}_s$ is 
\begin{align}
\left.
	 \begin{array}{l l} 
   \frac{\partial V}{\partial x}\!\! &= \frac{-F_2}{F_1-F_2}\frac{\partial t_{f_1}}{\partial x} + \frac{F_1}{F_1-F_2}\frac{\partial t_{f_2}}{\partial x}    \\
	&~~+ t_{f_1}\frac{\partial}{\partial x}\Big(\frac{-F_2}{F_1-F_2}\Big) + t_{f_2}\frac{\partial}{\partial x}\Big(\frac{F_1}{F_1-F_2}\Big) \\
	&~~+\frac{\partial V}{\partial \phi^*}\cdot \frac{d \phi^*}{d x}
\end{array}   \right.   \label{eq:GradRs}
\end{align}
where $t_{f_i}=t_{f_i}(\textbf{x})$ and $F_i=F_i(\textbf{x})$ are given by \eqref{eq:tfi} and \eqref{eq:Fi}, respectively, for $i=1,2$. 
The functions $t_{f_i}$, $F_i$, $\phi^*$, and $V$ depend only on the state $\textbf{x}$. The notation $(\textbf{x})$ is dropped hereafter in order to simplify the notation in \eqref{eq:GradRs} and in the remaining equations of this proof. In eq. \eqref{eq:GradRs} the variable $x$ denotes each element of the state $\textbf{x}$, that is, $x=\left\{x_E,y_E,x_{P_1},y_{P_1},x_{P_2},y_{P_2}\right\}$. The specific forms of $\frac{\partial t_{f_i}}{\partial x}$ are given by

\begin{align}
\left.
	 \begin{array}{l l} 
	 \frac{\partial t_{f_i}}{\partial x_E} &= \frac{1}{\beta_i^2-1}\Big[\cos\phi^* + \frac{1}{Q_i}\Big([(x_E-x_{P_i})\cos\phi^* \\
	&~~+ (y_E-y_{P_i})\sin\phi^*]\cos\phi^* \\
	&~~+(\beta_i^2-1)(x_E-x_{P_i})  \Big)   \Big]  	\\
	 \frac{\partial t_{f_i}}{\partial y_E} &= \frac{1}{\beta_i^2-1}\Big[\sin\phi^* + \frac{1}{Q_i}\Big([(x_E-x_{P_i})\cos\phi^*  \\
	&~~+ (y_E-y_{P_i})\sin\phi^*]\sin\phi^* \\
	&~~+(\beta_i^2-1)(y_E-y_{P_i})  \Big)   \Big]	 \\
	 \frac{\partial t_{f_i}}{\partial x_{P_j}} &=\left\{
	\begin{array}{l l}
	  -\frac{\partial t_{f_i}}{\partial x_E} \ \ \text{if} \ i=j\\
		0 \ \ \ \ \ \ \ \ \text{otherwise}
	\end{array}  \right. \\
 	 \frac{\partial t_{f_i}}{\partial y_{P_j}} &=\left\{
	\begin{array}{l l}
	  -\frac{\partial t_{f_i}}{\partial y_E} \ \ \text{if} \ i=j\\
		0 \ \ \ \ \ \ \ \ \text{otherwise}
	\end{array}  \right. 	 
\end{array}   \right.   \label{eq:Grads}
\end{align}
for $i,j=1,2$ where $Q_i=\big([(x_E-x_{P_i})\cos\phi^* + (y_E-y_{P_i})\sin\phi^*]^2 +(\beta_i^2-1)[(x_E-x_{P_i})^2 + (y_E-y_{P_i})^2] \big)^{1/2}$. Note that $Q_i>0$. 

Let us first analyze the terms in the second line of eq. \eqref{eq:GradRs}. Since $t_f=t_{f_1}=t_{f_2}$ we have that
 \begin{align}
\left.
	 \begin{array}{l l} 
t_{f_1}\frac{\partial}{\partial x}\Big(\frac{-F_2}{F_1-F_2}\Big) + t_{f_2}\frac{\partial}{\partial x}\Big(\frac{F_1}{F_1-F_2}\Big) 
	\!\!\!&= t_f \frac{\partial}{\partial x}\Big(\frac{F_1-F_2}{F_1-F_2}\Big)  \\
	&= t_f \cdot 0 \\ 
	&= \ 0.
\end{array}   \right.   \nonumber
\end{align}
We now evaluate the term $\frac{\partial V}{\partial \phi^*}$ as follows
 \begin{align}
\left.
	 \begin{array}{l l} 
\frac{\partial V}{\partial \phi^*}\!\!&=\frac{-F_2}{F_1-F_2}\frac{\partial t_{f_1}}{\partial \phi^*} + \frac{F_1}{F_1-F_2}\frac{\partial t_{f_2}}{\partial \phi^*}    \\ 
&~~+ t_{f_1}\frac{\partial}{\partial \phi^*}\big(\frac{-F_2}{F_1-F_2}\big) + t_{f_2}\frac{\partial}{\partial \phi^*}\big(\frac{F_1}{F_1-F_2}\big)  \\
&=\frac{-F_2}{F_1-F_2}\frac{\partial t_{f_1}}{\partial \phi^*} + \frac{F_1}{F_1-F_2}\frac{\partial t_{f_2}}{\partial \phi^*}    \\ 
&~~+ t_{f}\frac{\partial}{\partial \phi^*}\big(\frac{F_1-F_2}{F_1-F_2}\big)   \\
&=\frac{-F_2}{F_1-F_2}\frac{\partial t_{f_1}}{\partial \phi^*} + \frac{F_1}{F_1-F_2}\frac{\partial t_{f_2}}{\partial \phi^*}    \\ 
\end{array}   \right.   \nonumber
\end{align}
where $\frac{\partial t_{f_i}}{\partial \phi^*} = F_i \cdot t_{f_i}$. Therefore,
\begin{align}
\left.
	 \begin{array}{l l} 
\frac{\partial V}{\partial \phi^*}\!\!\!&=\frac{-F_1F_2}{F_1-F_2} t_{f_1} + \frac{F_1F_2}{F_1-F_2} t_{f_2}    
= \frac{F_1F_2 - F_1F_2}{F_1-F_2} t_f  
= 0
\end{array}   \right.   \nonumber
\end{align}
since $t_f= t_{f_1}= t_{f_2}$. Hence the gradient of the Value function \eqref{eq:GradRs} can be simplified and be written as
\begin{align}
\left.
	 \begin{array}{l l} 
   \frac{\partial V}{\partial x}\!\! &= \frac{-F_2}{F_1-F_2}\frac{\partial t_{f_1}}{\partial x} + \frac{F_1}{F_1-F_2}\frac{\partial t_{f_2}}{\partial x}    
\end{array}   \right.   \label{eq:GradRs2}
\end{align}
for $x=\left\{x_E,y_E,x_{P_1},y_{P_1},x_{P_2},y_{P_2}\right\}$. Note that $F_i(\textbf{x})$ can be written in the following form
\begin{align}
\left.
	 \begin{array}{l l} 
F_{i}(\textbf{x})\!\!&=\frac{(y_E-y_{P_i})\cos\phi^* - (x_E-x_{P_i})\sin\phi^* }{Q_i}  \\ 
&=\frac{\sin(\lambda_i-\phi^*)}{\sqrt{\cos^2(\lambda_i-\phi^*)+\beta_i^2 -1}}
\end{array}   \right.   \label{eq:Fi2}
\end{align}
for $i=1,2$. Since $\lambda_m \leq \phi^* \leq \lambda_M$, we have that $F_1$ and $F_2$ have different sign and $F_1-F_2\neq 0$. From \eqref{eq:Grads} and \eqref{eq:Fi2} we conclude that the Value function is $C^1$ inside region $\mathcal{R}_s$.

It is now critical to show that the Value function is $C^1$ on the boundary of each subset. We consider only the boundary between $\mathcal{R}_1$ and $\mathcal{R}_s$ since the proof for the boundary between $\mathcal{R}_2$ and $\mathcal{R}_s$ follows in the same way.
On the boundary between $\mathcal{R}_1$ and $\mathcal{R}_s$ we have that $\phi^*=\lambda_1$ because of \eqref{eq:Opt1}. We also have from \eqref{eq:Fi2} that $F_1(\phi^*=\lambda_1)=0$. Thus, the Value function evaluated in region $\mathcal{R}_s$ is given by
 \begin{align}
\left.
	 \begin{array}{l l} 
V(\textbf{x})&=\frac{-F_2}{-F_2} t_{f_1}   
= t_{f_1}.  
\end{array}   \right.   \nonumber
\end{align}
Then, it is sufficient to show that $t_{f_1}$ as given in $\mathcal{R}_s$ which is denoted by $t_{f_1s}$ in \eqref{eq:tfs1} below, is equal to $t_{f_1}$ as given in $\mathcal{R}_1$ (the subscript $s$ is attached to $t_{f_1}$ to show where $t_{f_1}$ is evaluated in $\mathcal{R}_s$). Consider \eqref{eq:tfi} evaluated at $\phi^*=\lambda_1$
 \begin{align}
\left.
	 \begin{array}{l l} 
t_{f_1s}\!\!\!&= \frac{1}{\beta_1^2-1} \Big[(x_E-x_{P1})\cos\lambda_1 + (y_E-y_{P1})\sin\lambda_1  \\
		&~~   + \big([(x_E-x_{P1})\cos\lambda_1 + (y_E-y_{P1})\sin\lambda_1]^2  \\
		&~~~~+(\beta_1^2-1)[(x_E-x_{P1})^2 + (y_E-y_{P1})^2] \big)^{1/2} \ \Big]  \\
     	  &= \frac{1}{\beta_1^2-1} \Big[\sqrt{(x_E-x_{P1})^2 + (y_E-y_{P1})^2}   \\
				&~~   + \big((x_E-x_{P1})^2 + (y_E-y_{P1})^2  \\
		&~~~~+(\beta_1^2-1)[(x_E-x_{P1})^2 + (y_E-y_{P1})^2] \big)^{1/2} \ \Big]  \\
	&=	\frac{(\beta_1+1)\sqrt{(x_E-x_{P1})^2+(y_E-y_{P1})^2}}{\beta^2_1-1}  \\
	&=\frac{\sqrt{(x_E-x_{P1})^2+(y_E-y_{P1})^2}}{\beta_1-1}	\\
	&=	t_{f_1}.
\end{array}   \right.  \label{eq:tfs1}
\end{align}
Thus, the Value function is continuous on the boundary between regions $\mathcal{R}_1$ and $\mathcal{R}_s$. 

Let us now consider $\frac{\partial V}{\partial \textbf{x}}$ on the same boundary. We want to show that \eqref{eq:GradRs2} is equal to \eqref{eq:Grad1} on the boundary between these two subsets, where the optimal Evader strategy is $\phi^*=\lambda_1$, and for the relevant states $x=\left\{x_E,y_E,x_{P1},y_{P1}\right\}$. Note that $\frac{\partial V}{\partial x_{P2}}=\frac{\partial V}{\partial y_{P2}}=0$ since $F_1=0$ and $\frac{\partial t_{f_1}}{\partial x_{P2}}=\frac{\partial t_{f_1}}{\partial y_{P2}}=0$. Let us consider first $\frac{\partial V}{\partial x_E}$ in \eqref{eq:GradRs2} evaluated at $\phi^*=\lambda_1$ 
\begin{align}
\left.
	 \begin{array}{l l} 
	 \frac{\partial V}{\partial x_E}\!\!\! &=  \frac{-F_2}{-F_2}\frac{\partial t_{f_1}}{\partial x_E} \\
	&=\frac{1}{\beta_1^2-1}\Big[\cos\lambda_1 \\
	&\ + {\scriptstyle {\frac{\sqrt{(x_E-x_{P1})^2 + (y_E-y_{P1})^2} \cos\lambda_1 +(\beta_1^2-1)(x_E-x_P)}{\beta_1\sqrt{(x_E-x_{P1})^2+(y_E-y_{P1})^2}}   \Big]} }%
	\\
	&= \frac{1}{\beta_1^2-1}\Big[\cos\lambda_1 + {\scriptstyle \frac{ \cos\lambda_1 +(\beta_1^2-1) \cos\lambda_1}{\beta_1} }%
		\Big]  \\
	&= 	\frac{\beta_1+1}{\beta_1^2-1}\cos\lambda_1  \\
	&= \frac{1}{\beta_1-1}\cos\lambda_1 
\end{array}   \right.   \label{eq:G1s}
\end{align}
which is equivalent to $\frac{\partial V}{\partial x_E}$ as given by \eqref{eq:Grad1}. Following a similar approach we obtain that $\frac{\partial V}{\partial y_E}$, $\frac{\partial V}{\partial x_{P1}}$, and $\frac{\partial V}{\partial y_{P1}}$ are also continuous on the boundary between regions $\mathcal{R}_1$ and $\mathcal{R}_s$. Hence the Value function is $C^1$ on the boundary between regions $\mathcal{R}_1$ and $\mathcal{R}_s$.

The same steps can be followed to show that the Value function is $C^1$ also on the boundary between regions $\mathcal{R}_2$ and $\mathcal{R}_s$. Therefore, altogether, the Value function is $C^1$ for $\textbf{x} \in\mathcal{R}$.

Finally, we show that the Value function satisfies the HJI equation 
\begin{equation}
\frac{\partial V({\bf x},t)}{\partial t} =
\arg \min_{\bf u_P} \max_{\bf u_E}
\left\{ \frac{\partial V({\bf x}, t)}{\partial {\bf x}} {\bf f(x, u_{\it P}, u_{\it E})} + 1 \right\}  \nonumber
\end{equation}
with ${\bf u}_{\it P} =[\psi_1, \psi_2]$ and ${\bf u}_{\it E}=\phi$. In this expression, the term 1 follows from the cost $\int_0^{t_f} 1 dt$. Then,
for the optimal feedback strategy, we have $\frac{\partial V({\bf x},t)}{\partial t}=0$ and $V({\bf x}, t)=V({\bf x})$; therefore, 
\begin{equation}
0 =\arg \min_{\bf u_P} \max_{\bf u_E}
\left\{ \frac{\partial V({\bf x})}{\partial {\bf x}} {\bf f(x, u_{\it P}, u_{\it E})} + 1 \right\}.
\end{equation}
If ${\bf u}^*_P={\bf u}^*_P({\bf x}) $ and ${\bf u}^*_E={\bf u}^*_E({\bf x})$ are (feedback) optimal solutions, then we can check if $V$ is the value function by checking   
\begin{equation}
0 =\left\{ \frac{\partial V({\bf x})}{\partial {\bf x}} {\bf f(x, u^*_{\it P}, u^*_{\it E})} + 1 \right\}.
\end{equation}
In Region 1 we have the following: $\phi^*=\psi_1^*=\lambda_1$ and
\begin{align}
\left.
	 \begin{array}{l l} 
	 \frac{\partial V}{\partial \textbf{x}} \cdot \textbf{f}(\textbf{x},\textbf{u}^*_P,\textbf{u}^*_E)  \\
	=\frac{\partial V}{\partial x_E}\cos\phi^* + \frac{\partial V}{\partial y_E}\sin\phi^*  \\
	~~+\beta_1\frac{\partial V}{\partial x_{P1}}\cos\psi_1^* +\beta_1\frac{\partial V}{\partial y_{P1}}\sin\psi_1^* \\
	 = {\scriptstyle \frac{1}{(\beta_1-1)\sqrt{(x_E-x_{P1})^2 + (y_E-y_{P1})^2}} }%
	  \big[(x_E-x_{P1})\cos\lambda_1  \\
	 ~~+ (y_E-y_{P1})\sin\lambda_1 -\beta_1(x_E-x_{P1})\cos\lambda_1 \\
	 ~~-\beta_1(y_E-y_{P1})\sin\lambda_1 \big]  \\
	 = {\scriptstyle \frac{1-\beta_1}{(\beta_1-1)\sqrt{(x_E-x_{P1})^2 + (y_E-y_{P1})^2}} }%
	\big[(x_E-x_{P1})\cos\lambda_1 \\
	~~+(y_E-y_{P1})\sin\lambda_1 \big] \\
   = {\scriptstyle \frac{-1}{\sqrt{(x_E-x_{P1})^2 + (y_E-y_{P1})^2}} }%
	\Big[ {\scriptstyle \frac{(x_E-x_{P1})^2 + (y_E-y_{P1})^2}{\sqrt{(x_E-x_{P1})^2 + (y_E-y_{P1})^2}} }%
	 \Big] \\
	= -1.
\end{array}   \right.   \label{eq:hji1}
\end{align}
Therefore, we have that
\begin{align}
\left.
	 \begin{array}{l l} 
	 1  + \frac{\partial V}{\partial \textbf{x}} \cdot \textbf{f}(\textbf{x},\textbf{u}^*_P,\textbf{u}^*_E)  
	= 1-1 = 0
\end{array}   \right.   \label{eq:hji1_2}
\end{align}
that is, the HJI equation is satisfied by the candidate Value function in Region $\mathcal{R}_1$. The same equation follows in $\mathcal{R}_2$. It is only left to show that the HJI equation is satisfied in Region $\mathcal{R}_s$. To accomplish this we use \eqref{eq:Grads} and \eqref{eq:GradRs2} to write 
\begin{align}
\left.
	 \begin{array}{l l} 
	& \frac{\partial V}{\partial \textbf{x}} \cdot \textbf{f}(\textbf{x},\textbf{u}^*_P,\textbf{u}^*_E)\\
	&= \frac{\partial V}{\partial x_E}\cos\phi^* + \frac{\partial V}{\partial y_E}\sin\phi^*  \\
	&~~+\beta_1\frac{\partial V}{\partial x_{P1}}\cos\psi_1^* +\beta_1\frac{\partial V}{\partial y_{P1}}\sin\psi_1^* \\
	&~~+\beta_2\frac{\partial V}{\partial x_{P2}}\cos\psi_2^* +\beta_2\frac{\partial V}{\partial y_{P2}}\sin\psi_2^* \\
	&=\Big(\frac{-F_2}{F_1-F_2}\frac{\partial t_{f_1}}{\partial x_E} + \frac{F_1}{F_1-F_2}\frac{\partial t_{f_2}}{\partial x_E}\Big) \cos\phi^* \\    
	&~~+\Big(\frac{-F_2}{F_1-F_2}\frac{\partial t_{f_1}}{\partial y_E} + \frac{F_1}{F_1-F_2}\frac{\partial t_{f_2}}{\partial y_E}\Big) \sin\phi^* \\    
	&~~+	\frac{-F_2}{F_1-F_2}\Big(\frac{\partial t_{f_1}}{\partial x_{P1}} \beta_1\cos\psi^*_1 + \frac{\partial t_{f_1}}{\partial y_{P1}}\beta_1\sin\psi^*_1\Big) \\  
		&~~+	\frac{F_1}{F_1-F_2}\Big(\frac{\partial t_{f_2}}{\partial x_{P2}} \beta_2\cos\psi^*_2 + \frac{\partial t_{f_2}}{\partial y_{P2}}\beta_2\sin\psi^*_2\Big) 
\end{array}   \right.   \label{eq:hjis}
\end{align}
subject to    
\begin{align}
\left.
	 \begin{array}{l l} 
	 \cos\psi_i^*&=\frac{\cos\phi^*}{\beta_i} + \frac{x_E-x_{P_i}}{\beta_i t_f}  \\
	 \sin\psi_i^*&=\frac{\sin\phi^*}{\beta_i} + \frac{y_E-y_{P_i}}{\beta_i t_f}
\end{array}   \right.   \label{eq:Consts}
\end{align}
for $i=1,2$. Note that \eqref{eq:Consts} are the optimal Pursuers' headings written in terms of the Evader's optimal heading $\phi^*$, where the Evader aims at the intersection point of the two Apollonius circles. In other words, \eqref{eq:Consts} are obtained from eqs. \eqref{eq:P1eqs} and \eqref{eq:P2eqs}.

Substituting \eqref{eq:Consts} into \eqref{eq:hjis} we obtain
\begin{align}
\left.
	 \begin{array}{l l} 
	 \frac{\partial V}{\partial \textbf{x}} \cdot \textbf{f}(\textbf{x},\textbf{u}^*_P,\textbf{u}^*_E)  \\
	=\frac{1}{t_f}\!\cdot\!\frac{-F_2}{F_1-F_2} \Big( \frac{\partial t_{f_1}}{\partial x_E}(x_{P1}\!-\!x_E) + \frac{\partial t_{f_1}}{\partial y_E}(y_{P1}\!-\!y_E) \Big) \\
	\ +\frac{1}{t_f}\!\cdot\!\frac{F_1}{F_1-F_2} \Big( \frac{\partial t_{f_2}}{\partial x_E}(x_{P2}\!-\!x_E) + \frac{\partial t_{f_2}}{\partial y_E}(y_{P2}\!-\!y_E) \Big) 
	\end{array}   \right.   \label{eq:HJBIs}
\end{align}
and we also evaluate the terms
	\begin{align}
\left.
	 \begin{array}{l l} 
	&\frac{1}{t_f} \Big( \frac{\partial t_{f_i}}{\partial x_E}(x_{P_i}-x_E) + \frac{\partial t_{f_i}}{\partial y_E}(y_{P_i}-y_E) \Big) \\
	&=-\frac{1}{t_f(\beta_i^2-1)}\Big[ (x_E-x_{P_i})\cos\phi^*+(y_E-y_{P_i})\sin\phi^*   \\
	&~~+\frac{1}{Q_i}\Big( [(x_E-x_{P_i})^2\cos^2\phi^* \\
	 &~~~~+ 2(x_E-x_{P_i})(y_E-y_{P_i})\sin\phi^*\cos\phi^* \\
	&~~~~+(y_E-y_{P_i})^2\sin^2\phi^*  \\
	&~~~~+ (\beta_i^2-1)[(x_E-x_{P_i})^2 + (y_E-y_{P_i})^2]\Big)\Big]  \\
	&= -\frac{1}{t_f(\beta_i^2-1)} \Big[ (x_E-x_{P_i})\cos\phi^*+(y_E-y_{P_i})\sin\phi^*   \\
	&~~ +\frac{1}{Q_i}\Big( [(x_E-x_{P_i})\cos\phi^* + (y_E-y_{P_i})\sin\phi^*]^2 \\
	&~~~~+(\beta_i^2-1)[(x_E-x_{P_i})^2 + (y_E-y_{P_i})^2]\Big)\Big] \\
	&= -\frac{1}{t_f(\beta_i^2-1)} \Big[ (x_E-x_{P_i})\cos\phi^*+(y_E-y_{P_i})\sin\phi^*   \\
	&~~~~+Q_i \Big] \\
	&= -\frac{1}{t_f} \cdot t_f  \\
	&=-1
\end{array}   \right.   \label{eq:hjis2}
\end{align}
for $i=1,2$. Substituting \eqref{eq:hjis2} into \eqref{eq:HJBIs} we have
\begin{align}
\left.
	 \begin{array}{l l} 
	 \frac{\partial V}{\partial \textbf{x}} \cdot \textbf{f}(\textbf{x},\textbf{u}^*_P,\textbf{u}^*_E)  
	&=\frac{-F_2}{F_1-F_2} \big( -1 \big) + \frac{F_1}{F_1-F_2} \big( -1 \big)   \\
	&= -1
	\end{array}   \right.  \nonumber   
\end{align}
and we can write
\begin{align}
\left.
	 \begin{array}{l l} 
	 1  + \frac{\partial V}{\partial \textbf{x}} \cdot \textbf{f}(\textbf{x},\textbf{u}^*_P,\textbf{u}^*_E)  
	= 1-1 = 0
\end{array}   \right.   \label{eq:hjis_2}
\end{align}
which means that the HJI equation is satisfied by the Value function in Region $\mathcal{R}_s$. $\square$

\textit{Remark}. The Value function \eqref{eq:ValueFn} determines the capture time when using the normalized (by $v_E$) speeds. In the case where $v_E\neq 1$ the Value of the game is simply scaled by $1/v_E$ to obtain the real capture time. We work with normalized speeds since the constant factor $1/v_E$ does not change the results in this paper.  

In summary, a solution of the 2P1EDG was proposed and the Value function was obtained. In this section we have verified that the Value function is continuous and continuously differentiable over $\mathcal{R}$. It was also shown that the candidate Value function satisfies the HJI equation for any $\textbf{x}\in\mathcal{R}$. All this together plus the uniqueness property of the Value function verifies that the proposed solution and the candidate Value function are in fact the solution of the Two-Pursuer One-Evader Differential Game.  


\subsection{Dispersal Surface and Multi-Pursuer Multi-Evader Differential Game}  \label{subsec:Dispersal}
The existence of a dispersal surface was discussed in Section \ref{subsec:Verification} and the condition is given by the distance from the Evader's position to each one of the two intersection of the Apollonius circles.
We illustrate the game's dispersal surface through an example.

Consider the initial positions:  $E=(5, 0)$,  $P_1=(0, 0)$,  $P_2=(24, -4)$. The players' speeds are: $v_E=1$, $v_{P_1}=1.25$, $v_{P_2}=1.3125$. 
Initially, there exist two sets of solutions that satisfy the optimality conditions; the two intersection points are at the same distance from the Evader's position and the state of the game is located on the singular surface of dispersal type, that is, $\textbf{x} \in\mathcal{D}$. 
In general, since there are more than one optimal solution and the each team does not know what strategy the opponent will implement, then two cases may occur: both teams choose the same solution or each team chooses a different solution. In the former case the state of the game stays on the dispersal surface and more than one solution is optimal. However, when teams choose different solutions then the state of the game leaves the dispersal surface and transitions into a regular region of the state space where only one optimal solution exists, hence, no dilemma exists any longer. Typically, Pursuers pay an infinitesimally small price for guessing incorrectly.

\begin{figure}
	\begin{center}
		\includegraphics[width=8.0cm,trim=.7cm .3cm .4cm .3cm]{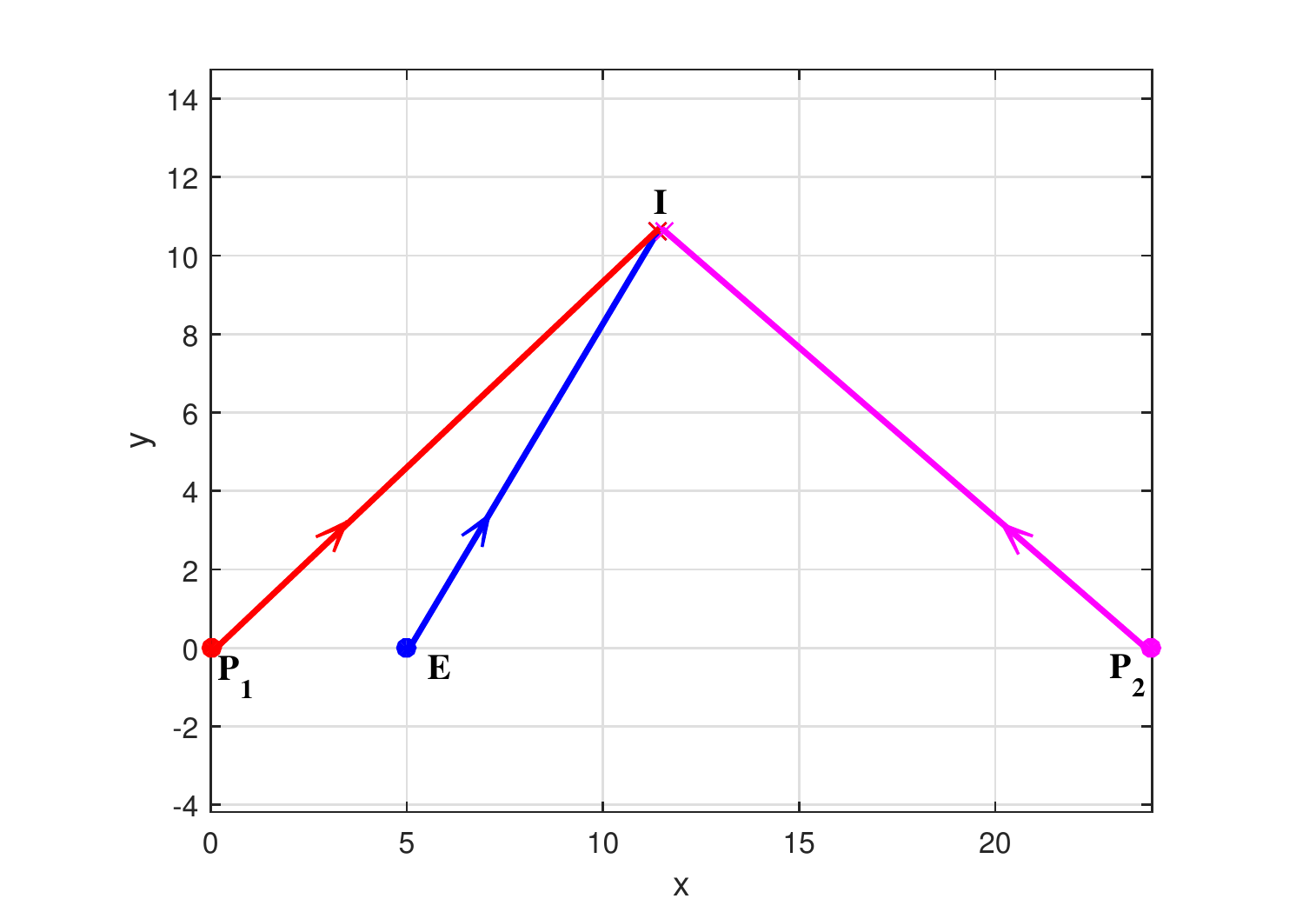}
	\caption{Optimal Trajectories}
	\label{fig:Ex2}
	\end{center}
\end{figure}

In this example each team makes a choice at the beginning of the game and it happens that they selected opposite solutions, the Pursuers initially head south whereas the Evader heads North. Because the players did not choose the same strategy, the state of the game leaves the dispersal surface and it transitions into the regular subspace $\mathcal{R}$. After measuring the Evader's new position, the Pursuers recompute the solution of the game; now it holds that $\textbf{x} \in \mathcal{R}_s^- \subset \mathcal{R}$.  Fig. \ref{fig:Ex2} shows the resulting trajectories. This problem illustrates a typical scenario in pursuit-evasion differential games where the pursuers generally see a small cost increment by incorrectly choosing between two or more strategies. This is inevitable in practice since the Evader will not announce its strategy and the Pursuers will pay a price for the state of the game to leave the dispersal surface and enter a regular subspace. The unlikely outcome where both teams choose the same strategy for all $0\leq t \leq t_f$ was referred by Isaacs as the perpetual dilemma \cite{Isaacs1965Differential}.

We have shown that the time of capture of a fleeing Evader can be reduced when two Pursuers work cooperatively as a team compared to the single Pursuer case. We can use the state feedback solution of the 2P1EDG in order to approximate the solution for a Multi-Pursuer Multi-Evader Differential Game, where, in addition of taking advantage of the cooperative behavior by the Pursuers, it is also necessary to assign teams of Pursuers to capture each Evader. The objective is to minimize the capture time of the last Evader where each Pursuer has to be assigned to one and only one Evader. Let us consider an example with five Pursuers and three Evaders where the speeds and initial positions of the five Pursuers are $v_P=\left\{1.3, 1.18, 1.2, 1.05, 1.1\right\}$ and $P=\left\{(3,9), (1,5), (0,0), (0.5,-3), (1.5,-7) \right\}$, respectively. The initial positions and speeds of the three Evaders are $v_E=\left\{0.98, 0.85, 0.76\right\}$ and $E=\left\{(8,5), (10, 1), (7,-3) \right\}$, respectively. 

In this example we consider assignments of teams $2-2-1$. We do not consider cases where three Pursuers may be assigned to a single Evader. Under these conditions we obtain the combination of possible assignments shown in Table \ref{table:Com}, where each numerical entry denotes the Value of the Game for the 2P1EDG combination $P_iP_jE_l$ and the superscript denotes which case is active under such combination, that is, if the Evader will be captured only by $P_i$, only by $P_j$, or simultaneously by both Pursuers. A given Pursuer can be assigned to capture only one Evader which imposes constraints on the selection of assignments; if for instance, $P_1$ is assigned to $E_1$, either alone or teaming up with $P_5$, then it cannot be assigned to $E_2$ (or $E_3$) and the best time to capture $E_2$ would be 28.46. Similarly, if $P_1$ is assigned either to $E_2$ or $E_3$, then it cannot be used to capture $E_1$ and the best time to capture $E_1$ would be 35.00.

 The best assignment strategy is $\left\{P_1E_1\right\}$, $\left\{P_2P_3E_2\right\}$, $\left\{P_4P_5E_3\right\}$ corresponding to the capture times 20.01, 28.46, and 19.97, respectively. All Evaders have been captured at 28.46 time units and no other combination can reduce this value.

\begin{table}
	\caption{Possible combinations of assignments $P_i,P_j,E_l$ for $i,j=1,...,5$ and $l=1,2,3$. The superscript denotes the active case}
\begin{center}
  \begin{tabular}{| c | c | c | c |}
    \hline
     {$P_i,P_j$} \ - \ {$E_l$} & 1 & 2 & 3 \\ \hline
      1,2 & $20.01^{(1)}$ & $23.62^{(1)}$ & $23.36^{(s)}$ \\  \hline
			1,3 & $20.01^{(1)}$ & $23.33^{(s)}$ & $17.31^{(3)}$  \\  \hline
			1,4 & $20.01^{(1)}$ & $23.62^{(1)}$ & $20.18^{(s)}$ \\  \hline
			1,5 & $20.00^{(s)}$ & $22.92^{(s)}$ & $16.34^{(s)}$ \\  \hline
      2,3 & $35.00^{(2)}$ & $28.46^{(s)}$ & $17.31^{(3)}$ \\  \hline
      2,4 & $35.00^{(2)}$ & $29.84^{(2)}$ & $21.41^{(s)}$  \\  \hline  
			2,5 & $35.00^{(2)}$ & $29.68^{(s)}$ & $17.66^{(s)}$ \\  \hline
			3,4 & $42.88^{(3)}$ & $28.71^{(3)}$ & $17.31^{(3)}$ \\  \hline
      3,5 & $42.88^{(3)}$ & $28.71^{(3)}$ & $16.75^{(s)}$ \\  \hline
      4,5 & $113.73^{(5)}$ & $46.69^{(5)}$ & $19.97^{(s)}$  \\  \hline  
  \end{tabular} 	
	\label{table:Com}
\end{center}
\end{table}

\section{Active Target Defense Differential Games (ATDDG)} \label{sec:litATDDG}
\subsection{Overview}
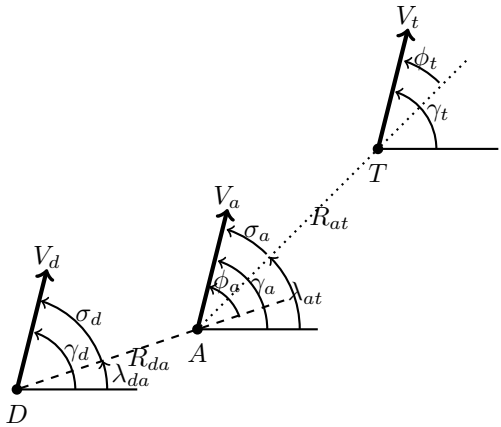
\begin{figure}[htb]
\begin{center}
\begin{tikzpicture}[scale = 0.8,thick]
\node (defender) 	at (0,0) [circ] {};
\node (attacker) 	at (3,1) [circ] {};
\node (target)		at (6,4) [circ] {};
\draw (defender) node [below = 2pt] {$D$};
\draw (attacker) node [below = 2pt] {$A$};
\draw (target)   node [below = 2pt] {$T$};
\draw[solid] (defender) -- ($(defender) + (2,0)$);
\draw[solid] (attacker) -- ($(attacker) + (2,0)$);
\draw[solid] (target) -- ($(target) + (2,0)$);
\draw[dashed] (defender)--(attacker);
\draw[dashed] (attacker)--($(attacker)+0.5*(3,1)$);
\draw[dotted] (attacker)--(target);
\draw[dotted] (target)--($(target) + 0.5*(3,3)$);
\draw[solid,ultra thick,->] (defender)--($(defender) + (0.5,2.0)$);
\draw[solid,ultra thick,->] (attacker)--($(attacker) + (0.5,2.0)$);
\draw[solid,ultra thick,->] (target) -- ($(target) + (0.5,2.0)$);
\draw [solid,thick,->] ($(defender)+(0.97,0)$) arc (0:71.56:1);
\draw [solid,thick,->] ($(attacker)+(1.16,0)$) arc (0:71.56:1.2);
\draw [solid,thick,->] ($(target)+(0.97,0)$)   arc (0:71.56:1);
\draw [solid,thick,->] ($(defender)+(1.5,0)$) arc (0:18.43:1.5);
\draw [solid,thick,->] ($(attacker)+(1.7,0)$) arc (0:45:1.7);
\draw [solid,thick,->] ($(defender)+(1.41,0.46)$) arc (18.43:71.56:1.6);
\draw [solid,thick,->] ($(attacker)+(1.15,1.25)$) arc (45:71.56:1.7);
\draw [solid,thick,->] ($(target)+(1.02,1.1)$) arc (45:71.56:1.5);
\draw [solid,thick,->] ($(attacker)+(0.7,0.2)$) arc (18.43:71.56:0.8);
\draw ($(defender) + (1.9,0.25)$)  node [] {$\lambda_{da}$};
\draw ($(attacker) + (1.8,0.6)$)   node [] {$\lambda_{at}$};
\draw ($(defender) + (1.0,0.6)$)   node [] {$\gamma_{d}$};
\draw ($(attacker) + (1.1,0.7)$)   node [] {$\gamma_{a}$};
\draw ($(target)   + (1.0,0.6)$)   node [] {$\gamma_{t}$};
\draw ($(defender) + (1.2,1.2)$)   node [] {$\sigma_{d}$};
\draw ($(attacker) + (1.0,1.6)$)   node [] {$\sigma_{a}$};
\draw ($(attacker) + (0.5,0.8)$)   node [] {$\phi_{a}$};
\draw ($(target)   + (0.8,1.5)$)   node [] {$\phi_{t}$};
\draw ($(defender) + (2.2,0.5)$)   node [] {$R_{da}$};
\draw ($(attacker) + (2.2,1.85)$)   node [] {$R_{at}$};
\draw ($(defender) + (0.5,2.2)$)   node [] {$V_{d}$};
\draw ($(attacker) + (0.5,2.2)$)   node [] {$V_{a}$};
\draw ($(target) + (0.5,2.2)$)   node [] {$V_{t}$};
\end{tikzpicture}
\caption{Defender-Attacker-Target Geometry describing the Active Target Defense Differential Games where the Defender pursues the Attacker who pursues the actively maneuvering Target, \cite{Weintraub2017Optimal, Weintraub2018Optimal}.}
\label{fig:DAT}
\end{center}
\end{figure}

Target Defense Differential Games (TDDG) are recently introduced pursuit-evasion differential games with three agents. A target ($T$) is pursued by an agent called the attacker ($A$). A third player, the Defender ($D$) pursues $A$ in an effort to defend $T$. The outcome of the three-player game is simple: If $D$ intercepts $A$ before $A$ captures $T$, then the Target is successfully defended; however, if $A$ captures $T$ before $D$ can intercept $A$, then the defense is unsuccessful and $A$ is the winner. When the target is able to maneuver in the three-player game, then $T$ and $D$ cooperate while playing against $A$ and this scenario is known as the Active Target Defense Differential Game (ATDDG). Fig.~\ref{fig:DAT}
 describes the geometry used to describe ATDDGs \cite{Weintraub2017Optimal, Weintraub2018Optimal}.
 
There exist a number of popular performance metrics which are used when posing the ATDDG. In games of kind, the interest lies with the outcome of the defense: does the Target succeed in evading the attacker or is captured? In games of degree, when the target escapes, the range from the attacker to the target at the instant when the attacker is intercepted by the defender has been used. When the target is captured by the Attacker, the distance between the defender and attacker at the instant of capture of the Target by the Attacker has been used as a performance metric. These ranging metrics are popular because they quantify the outcome of the engagement and are readily computed. Other metrics, such as time to capture, are also of interest.

The two-players one target game of \cite{Cardaliaguett1996a} is an early version of target defense. In the paper, a two-player differential game is played wherein one of the players wants the state of the system to reach a target, while the other player wants the state of the system to avoid this target. Introduced by Boyell, the defense of a ship from an incoming torpedo using a counter-weapon was described, \cite{Boyell1976,Boyell1980Counterweapon}. The authors of \cite{Yamasaki2010, Yamasaki2013} proposed the defense of an active target by launching a defensive missile. 
Rubinsky and Gutman presented an analysis of the end-game ATDDG scenario based on the Attacker-Target miss distance for a non-cooperative Target-Defender. The authors develop linearization-based Attacker maneuvers in order to evade the Defender and continue pursuing the Target using the LQ paradigm, \cite{Rubinsky2012, Rubinsky2014}. Defense of  non-maneuverable aircraft was addressed in \cite{Venkatesan2015New,Weintraub2018Kinematic}.
In \cite{Rusnak2011}, the limiting values of the three participants optimal strategies are studied as the quadratic weight on the defending missiles acceleration command tends to zero. They show that in the limit, the intercepting missiles and the targets optimal strategies are identical in form to that obtained in the game without the defending missile. Ratnoo and Shima proposed a game theoretic analysis of the ATDDG problem using conventional guidance laws for both attacker and defender, \cite{Ratnoo2011, Ratnoo2012}. The cooperative strategies proposed by \cite{Shima2011} allowed for a maneuverability disadvantage of the Defender with respect to the Attacker and the results show that the optimal target maneuver is either constant or arbitrary. The authors of \cite{Shaferman2010} implemented a Multiple Model Adaptive Estimator (MMAE) to identify the guidance law and parameters of incoming missiles and optimize defender strategies to minimize the control effort. 
In \cite{Sun2017}  the three-agent game using zero-effort-miss (ZEM) is investigated.

The ATDDG has been addressed in great detail in
\cite{Pachter2014Active, Garcia2015Active, Casbeer2018, Weintraub2017Optimal, Weintraub2018Optimal, Garcia2017e, Garcia2017f}. The ATDDG was analyzed  in \cite{Pachter2014Active, Garcia2015Active} using simple motion kinematics. In order to find the optimal strategy for the Defender to intercept the Attacker, the geometric concept of the Apollonius circle is used. In this analysis, the authors are able to look at the critical Target/Attacker speed ratio to ensure target survival. 

The authors of \cite{Casbeer2018} considered the use of two defenders to better engage the attacker. The work in  \cite{Garcia2017f} investigated the case where the Defender had a non-zero capture radius. While in previous work, information between all players was shared, the work in \cite{Weintraub2017Optimal, Garcia2017e} was concerned with optimal evasion of the target assuming the defender's and attacker's control laws were proportional navigation and pure pursuit, respectively. Information was restricted to the Attacker and Defender, and heading rate constraints were imposed on the Target. Finally, in \cite{Weintraub2018Optimal}, an Extended Kalman Filter was used to investigate the same engagement with sensor models.

\subsection{Problem Formulation}
The scenario of active target defense considers three players: a Target ($T$), an Attacker ($A$), and a Defender ($D$) which have ``simple motion" $\grave{\text{a}}$ la Isaacs.  The game is played in the Euclidean plane where the controls of $T$, $A$, and $D$ are their respective instantaneous headings $\phi$, $\chi$, and $\psi$ and their states are specified by their Cartesian coordinates $\textbf{x}_T=(x_T,y_T)$, $\textbf{x}_A=(x_A,y_A)$, and $\textbf{x}_D=(x_D,y_D)$. The players $T$, $A$ and $D$ have constant speeds denoted by $V_T$, $V_A$, and $V_D$, respectively.  The complete state of the TAD differential game is specified by $\textbf{x}:=( x_T, y_T, x_A, y_A, x_D, y_D)\in \mathbb{R}^6$. The game set is the entire space $\mathbb{R}^6$. The initial time is denoted by $t_0$ and the corresponding initial state of the system is $\textbf{x}_0 := (x_{T_0}, y_{T_0}, x_{A_0}, y_{A_0}, x_{D_0}, y_{D_0}) = \textbf{x}(t_0)$.
The Target aircraft is slower than the Attacker missile, and thus the speed ratio $\alpha=V_T/V_A<1$. We assume that the Attacker and Defender have similar capabilities, so $V_A=V_D$. Without loss of generality, the players' speeds are normalized so that $V_A=V_D=1$ and $V_T=\alpha$.

 The control input of the $T/D$ team is the pair of instantaneous headings $\textbf{u}_{T,D}=\left\{\phi,\psi\right\}$.
The Attacker's control is his instantaneous heading angle, $\textbf{u}_A=\left\{\chi\right\}$.  The dynamics $\dot{\textbf{x}}=\textbf{f}(\textbf{x},\textbf{u}_A,\textbf{u}_{T,D})$ are specified by the system of ordinary differential equations
\begin{align}
 \left.
	\begin{array}{l l}
        \dot{x}_A&=\cos\chi,    \ \ \ \ \ x_A(0)=x_{A_0}   \\
	\dot{y}_A&=\sin\chi,  \ \ \ \ \ y_A(0)=y_{A_0}  \\
	\dot{x}_D&=\cos\psi,   \ \ \ \ \ x_D(0)=x_{D_0} \\
	\dot{y}_D&=\sin\psi,   \ \ \ \ \ y_D(0)=y_{D_0}  \\
	\dot{x}_T&=\alpha\cos\phi,  \ \ \  x_T(0)=x_{T_0}      \\
	\dot{y}_T&=\alpha\sin\phi,    \ \ \  y_T(0)=y_{T_0}    
	\end{array}  \right.   \label{eq:xA}
\end{align}
where the admissible controls are given by $\chi,\phi,\psi \in [-\pi,\pi]$. Both, the state and the controls, are unconstrained.
The terminal condition is point capture, that is, the separation between Target and Attacker becomes zero allowing the Attacker to capture the Target and win the game. An alternative termination condition is  when the separation between Attacker and Defender is equal to zero; this case represents interception of the Attacker by the Defender  and the $T/D$ team wins. Hence, the termination set for the complete TAD differential game is
\begin{align}
   \mathcal{S} :=  \mathcal{S}_e    \   \bigcup \   \mathcal{S}_c   \label{eq:TwoSets}
\end{align}
where 
\begin{align}
   \mathcal{S}_e:=  \big\{ \ \textbf{x} \ | \sqrt{ (x_A-x_D)^2 + (y_A-y_D)^2} =0  \big\}   \label{eq:escape}
\end{align}
represents interception of the Attacker by the Defender (and the Target escapes) and 
\begin{align}
   \mathcal{S}_c:=   \big\{ \ \textbf{x} \ |  \sqrt{ (x_A-x_T)^2 + (y_A-y_T)^2}=0 \big\}    \label{eq:capture}
\end{align}
represents the opposite outcome where the Attacker wins by capturing the Target. 

In this paper the focus is on the Game of Degree when \eqref{eq:escape} applies. In such a case the cost functional is given by
\begin{align}
  J_c(\textbf{u}_A(t),\textbf{u}_B(t),\tilde{\textbf{x}}_0):=\frac{1}{2}[(x_{A_f}-x_{T_f})^2+(y_{A_f}-y_{T_f})^2]   \label{eq:EscCost}
\end{align}
subject to the terminal condition \eqref{eq:escape}.

The next two sections summarize the approach described in  \cite{Pachter19} where a reduced state space is used to derive the optimal strategies of the escape game. Without loss of generalization, it is assumed that $x_{D}=-x_{A}$, $y_{A}=0$, $y_{D}=0$, and $y_T\geq 0$; the scenario is as shown in Fig. \ref{fig:PMP}. The selected reduced state space is a suitable choice of framework for this particular problem. Other choices such as a framework based on relative distances is common in other games. In general, the main purpose is to simplify the analysis by reducing the dimension of the system.

\subsection{Game of Kind}

The concept of Game of Kind is fundamental in differential games and, given the initial conditions and problem parameters, its solution provides the answer to the question of which team wins the game.

In the ATDDG and employing the reduced state space $\textbf{x}=(x_A,x_T,y_T)$ we have that when $x_T<0$ the Target can escape. When $x_T>0$ the state space is partitioned into two regions: $R_e$ and $R_{c}$. The region $R_e$ is the set of states such that if the Target's initial position $(x_T,y_T)$ and the coordinate $x_A$ are such that $(x_A,x_T,y_T)\in R_e$, then, the Target is guaranteed to escape the Attacker, provided the Target and the Defender team implement their optimal strategies $\phi^*$ and $\psi^*$.  The Game of Degree, that is, the ATDDG, is played in $R_e$. The region $R_{c}$ specifies the states where under optimal Attacker play, notwithstanding the presence of the Defender, the Target cannot escape.

For a specified speed ratio $0<\alpha<1$, the equation
\begin{align}
     \alpha=\frac{\sqrt{(x_A+x_T)^2+y_T^2}-\sqrt{(x_A-x_T)^2+y_T^2}}{2x_A}    \label{eq:alphasol1}
\end{align}
renders the boundary in the state space $(x_A,x_T,y_T)$, $x_T,x_A>0$,  where the active target defense game of degree is played out, and as such this boundary constitutes the solution to the Game of Kind.

\begin{proposition}
 \cite{Pachter19}. For a given speed ratio $0<\alpha<1$ the region of win of the Attacker is
\begin{align}
R_c = \{ (x_A,x_T,y_T)| x_A>0, \ x_T>0, \ y_T\geq 0,  \nonumber  \\
  \ x_A^2 + \frac{y_T^2}{1-\alpha^2} - \frac{x_T^2}{\alpha^2}<0 \}.  \nonumber
\end{align}
The manifold $\mathcal{B}$ in $\mathbb{R}^3$ which is the boundary of the set $R_e$ for which the Target is guaranteed to escape is 
\begin{align}
\mathcal{B} = \{ (x_A,x_T,y_T)| x_A>0, \ x_T>0, \ y_T\geq 0,  \nonumber  \\
  \ x_A^2 + \frac{y_T^2}{1-\alpha^2} - \frac{x_T^2}{\alpha^2}=0 \}.     \nonumber
\end{align}
and the region of win of the $T$ \& $D$ team where the ATDDG is played is  
\begin{align}
R_e = \{ (x_A,x_T,y_T)| x_A\geq 0, \ x_T \geq 0, \ y_T \geq 0,   \nonumber   \\
  \ x_A^2 \!+\! \frac{y_T^2}{1\!-\!\alpha^2} \!-\! \frac{x_T^2}{\alpha^2}>0 \}  \cup  \{ (x_A,x_T,y_T)| x_A\geq 0, \ x_T \leq 0 \}.  \nonumber
\end{align}
Thus, for a fixed $x_A>0$, the $x_A$-cross section of $\mathcal{B}$ which divides the reduced state space into the two regions $R_e$ and $R_{c}$ is the right branch of the hyperbola (where, $x_T>0$)
\begin{align}
	  \frac{x_T^2}{\alpha^2x_A^2} -\frac{y_T^2}{(1-\alpha^2)x_A^2} = 1.     \label{eq:hb}
\end{align}
\end{proposition}

\subsection{Game of Degree}
In the following Theorem, the ``two-sided''  PMP is used in order to synthesize the players' optimal strategies. This problem illustrates how the necessary conditions for optimality or PMP can be used to obtain the regular solution of this differential game where $T$ and $D$ cooperate to maximize the terminal range while $A$ aims at minimizing it.

\begin{theorem}  \label{th:main}
Consider the ATDDG \eqref{eq:xA}-\eqref{eq:EscCost}. The problem parameter is the speed ratio $0\leq\alpha < 1$ and the reduced state space, is $(x_A,x_T,y_T)$ where, without loss of generality, $x_A>0$ and $y_T\geq 0$. In the part of the state space where the Target can escape, the optimal headings of the Attacker, the Target, and the Defender are given by the state feedback control laws
\begin{align}
 \left.
	 \begin{array}{l l}
\cos\phi^*=\pm \frac{x_{T}}{\sqrt{x_{T}^2+(y_{T}-y)^2}}  \\
 \sin\phi^*=\pm \frac{y_{T}-y}{\sqrt{x_{T}^2+(y_{T}-y)^2}} \ \  \text{for} \ x_T\neq 0.     \\ 
  \cos\chi^*=-\frac{x_{A}}{\sqrt{x_{A}^2+y^2}}, \ \ \   \sin\chi^*=\frac{y}{\sqrt{x_{A}^2+y^2}}.      \\
 \cos\psi^*=\frac{x_{A}}{\sqrt{x_{A}^2+y^2}}, \ \ \ \    \sin\psi^*=\frac{y}{\sqrt{x_{A}^2+y^2}}.    
\end{array}  \right.  \label{eq:OC}
\end{align}
where $y$ is a real solution of the quartic equation
\begin{align}
   \left.
	 \begin{array}{l l}
     (1\!-\!\alpha^2)y^4 - 2(1\!-\!\alpha^2)y_Ty^3   + \\
		 \big[(1\!-\!\alpha^2)y_T^2\!+\!x_A^2\!-\!\alpha^2x_T^2\big]y^2  \!-\! 2x_A^2y_Ty\!+\!x_A^2y_T^2\!=\!0
	\end{array}  \label{eq:Quartic}  \right.
\end{align}
which is parameterized by the speed ratio $0\leq\alpha<1$. The quartic equation has two real solutions $y_1$ and $y_2$ and $y_1\leq y_T\leq y_2$. In eqs. \eqref{eq:OC} when $x_T\leq 0$ the solution $y_1$ is selected and when $x_T>0$ the solution $y_2$ is used. Also, when $x_T<0$, the sign in the Target heading in eq. \eqref{eq:OC} is positive and, when $x_T>0$, the sign in the Target heading is negative. Finally,  when $x_T=0$, the Target's optimal heading, $\phi^*$,  is given by 
\begin{align}
  \phi^*=\varphi^*+\pi/2
\end{align}
where
\begin{align}
  \tan \varphi^* = \frac{\sqrt{x_A^2+(1-\alpha^2)y_T^2}}{\alpha y_T}.  	\label{eq:tanvarphi}
\end{align}
\end{theorem}
\textit{Proof}. The optimal control inputs in terms of the co-state variables are obtained from Isaacs' Main Equation 1 $ \min_{\phi,\psi} \max_\chi \mathcal{H} =0$ and they are characterized by
\begin{align}
  \left.
	 \begin{array}{l l}
   \cos\chi^*=\frac{\lambda_{x_A}}{\sqrt{\lambda_{x_A}^2+\lambda_{y_A}^2}}, \ \ \ \ \  \sin\chi^*=\frac{\lambda_{y_A}}{\sqrt{\lambda_{x_A}^2+\lambda_{y_A}^2}}  \nonumber  \\ 
  \cos\psi^*=-\frac{\lambda_{x_D}}{\sqrt{\lambda_{x_D}^2+\lambda_{y_D}^2}}, \ \   \sin\psi^*=-\frac{\lambda_{y_D}}{\sqrt{\lambda_{x_D}^2+\lambda_{y_D}^2}} \nonumber \\ 
	  \cos\phi^*=-\frac{\lambda_{x_T}}{\sqrt{\lambda_{x_T}^2+\lambda_{y_T}^2}}, \ \ \   \sin\phi^*=-\frac{\lambda_{y_T}}{\sqrt{\lambda_{x_T}^2+\lambda_{y_T}^2}}.  \nonumber 
\end{array}  \right.
\end{align}
Additionally, the co-state dynamics are $\dot{\lambda}_{x_A}=\dot{\lambda}_{y_A}=\dot{\lambda}_{x_D}=\dot{\lambda}_{y_D}=\dot{\lambda}_{x_T}=\dot{\lambda}_{y_T}=0$.
 Hence, all co-states are constant and we have that $\chi^*\equiv$ constant, $\psi^*\equiv$ constant, and $\phi^*\equiv$ constant. In other words, the regular optimal trajectories are straight lines.

Concerning the solution of the attendant Two-Point Boundary Value Problem (TPBVP) on $0\leq t\leq t_f$ in $\mathbb{R}^{12}$, we have $6$ initial states specified by \eqref{eq:xA} and we need $6$ more conditions for the terminal time $t_f$.
In this respect, define the augmented Mayer terminal value function $\Phi_a:\mathbb{R}^6\rightarrow \mathbb{R}^1$ 
\begin{align}
 \left.
	 \begin{array}{l l}
\Phi_a(\textbf{x}_f):= \frac{1}{2}[(x_A(t_f)-x_T(t_f))^2+(y_A(t_f)-y_T(t_f))^2]  \\
\qquad \qquad  \  +\nu_1(x_A(t_f)-x_D(t_f))+\nu_2(y_A(t_f)-y_D(t_f))
\end{array}  \right.  \nonumber
\end{align}
where $\nu_1$ and $\nu_2$ are Lagrange multipliers. The PMP or Dynamic Programming yields the transversality/terminal conditions
$\lambda(t_f)=-\frac{\partial}{\partial \textbf{x}}\Phi_a(\textbf{x}_f)$, that is,
\begin{align}
&\lambda_{x_A}=x_T(t_f)-x_A(t_f)-\nu_1    \label{eq:Lxa}   \\
&\lambda_{y_A}=y_T(t_f)-y_A(t_f)-\nu_2    \label{eq:Lya}   \\
&\lambda_{x_D}=\nu_1    \label{eq:Lxd}   \\
&\lambda_{y_D}=\nu_2    \label{eq:Lyd}   \\
&\lambda_{x_T}=x_A(t_f)-x_T(t_f)    \label{eq:Lxt}   \\
&\lambda_{y_T}=y_A(t_f)-y_T(t_f)    \label{eq:Lyt}   
\end{align}
At this point, we have that equations \eqref{eq:Lxa}-\eqref{eq:Lyt} together with the terminal equations
\begin{align}
   x_A=x_D, \ \ \   y_A=y_D.   \label{eq:PCx} 
\end{align}
 yield $8$ conditions. Since we need only $6$ conditions we eliminate the introduced Lagrange multipliers $\nu_1$ and $\nu_2$ from equations \eqref{eq:Lxa}-\eqref{eq:Lyt} and we obtain 
\begin{align}
&\lambda_{x_A}+\lambda_{x_D}=x_T(t_f)-x_A(t_f)    \label{eq:Lxad}   \\
&\lambda_{y_A}+\lambda_{y_D}=y_T(t_f)-y_A(t_f)    \label{eq:Lyad}   
\end{align}
Finally, the time $t_f$ is specified by the PMP requirement that the Hamiltonian $\mathcal{H}(\textbf{x}(t),\lambda(t),\chi,\psi,\phi) |_{t_f}\equiv 0$. 


\begin{figure}
	\begin{center}
		\includegraphics[width=8.4cm,height=4.0cm,trim=.5cm 2.8cm .2cm .8cm]{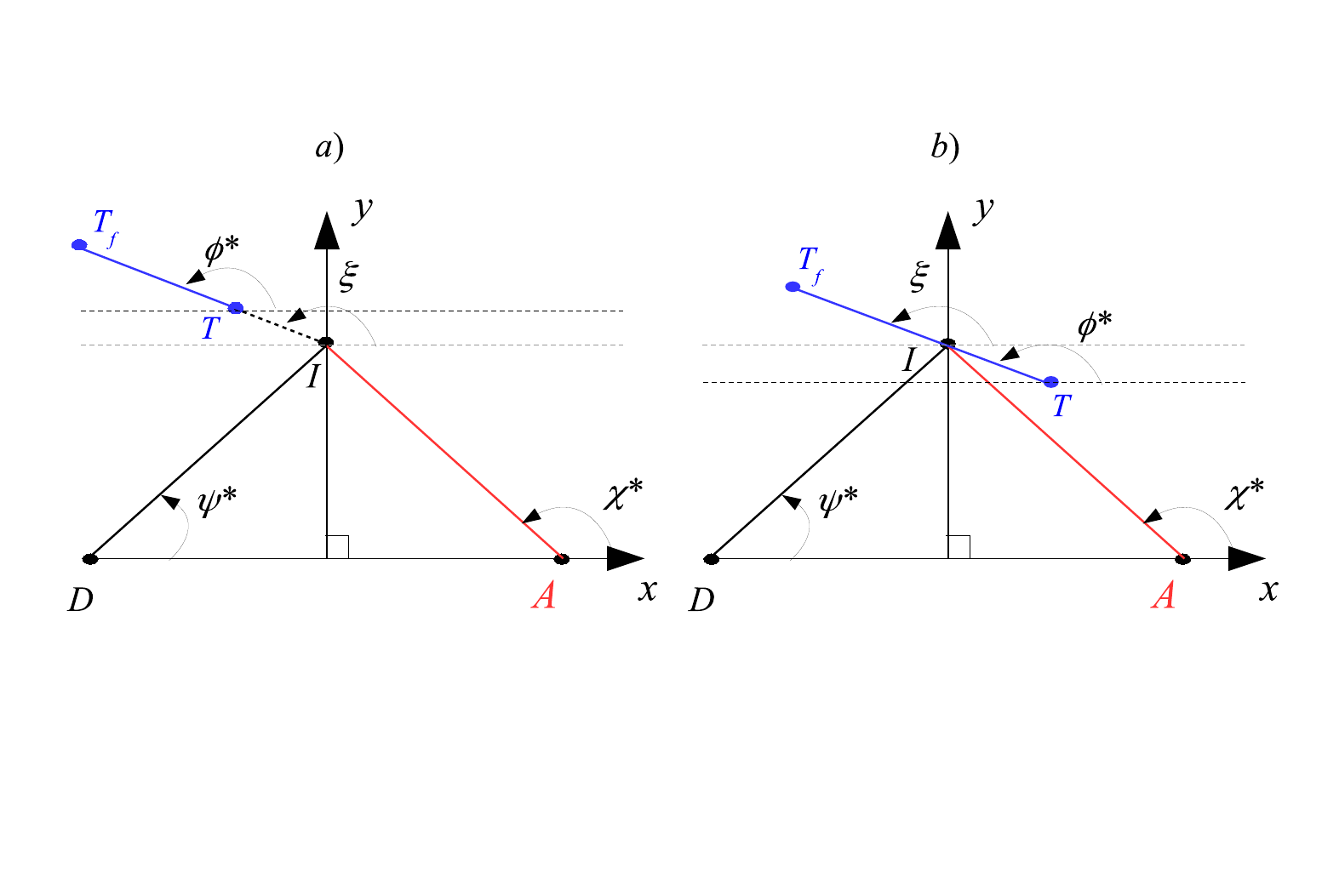}
	\caption{Optimal Headings of the Target, the Attacker, and the Defender}
	\label{fig:PMP}
	\end{center}
\end{figure}

Because the optimal trajectories of $A$, $D$, and $T$ are straight lines and $V_D=V_A$ we have that 
\begin{align}
&x_A(t_f)=0      \\
&x_D(t_f)=0     \\
&y_A(t_f)=y_D(t_f)
\end{align}
Let $y:= y_A(t_f)=y_D(t_f)$. Also, let $x_A=x_A(t')$, $x_T=x_T(t')$, and $y_T=y_T(t')$ be the instantaneous positions at some time $t'<t_f$. Hence, from equations \eqref{eq:xA} we obtain the following
\begin{align}
	&x_T(t_f)=x_{T}+\alpha \cdot (t_f-t')\cos\phi,    \label{eq:xTt}  \\
	&y_T(t_f)=y_{T}+\alpha \cdot (t_f-t')\sin\phi,     \label{eq:yTt}  \\
  &0=x_{A}+(t_f-t')\cos\chi,      \label{eq:xAt}  \\
	&y=(t_f-t')\sin\chi,   \label{eq:yAt} \\
	&0=-x_{A}+(t_f-t')\cos\psi,    \label{eq:xDt}  \\
	&y=(t_f-t')\sin\psi,    	\label{eq:yDt}
\end{align}
In addition, equations \eqref{eq:Lxt}-\eqref{eq:Lyad} can be written as follows
\begin{align}
&\lambda_{x_T}=-x_T(t_f)  \label{eq:Lxtt} \\
&\lambda_{y_T}=y-y_T(t_f)  \label{eq:Lytt} \\
&\lambda_{x_A}+\lambda_{x_D}=x_T(t_f)   \label{eq:Lxadt}   \\
&\lambda_{y_A}+\lambda_{y_D}=y_T(t_f)-y    \label{eq:Lyadt}   
\end{align}
Consider eq. \eqref{eq:phico} and eqs. \eqref{eq:Lxtt} and \eqref{eq:Lytt}. We calculate the following 
\begin{align}
 \left.
	 \begin{array}{l l}
\cos\phi^*=\frac{x_T(t_f)}{\sqrt{x_T^2(t_f)+(y-y_T(t_f))^2}} \\
  \sin\phi^*=\frac{y_T(t_f)-y}{\sqrt{x_T^2(t_f)+(y-y_T(t_f))^2}}.  \label{eq:phigeo} 
\end{array}  \right.
\end{align}
Hence, in light of eqs. \eqref{eq:xTt} and \eqref{eq:yTt}, we have that the optimal heading angle of $T$ is $ \phi^*=\xi $
where the angle $\xi$ is shown in Fig. \ref{fig:PMP}. The three points, $T$, $T_f$, and $y$ are collinear and the optimal geometry is as shown in Fig. \ref{fig:PMP} where $T_f=(x_T(t_f),y_T(t_f))$.
From the $\triangle ADy$ in Fig. \ref{fig:PMP} we conclude that $t_f-t'=\sqrt{x_{A}^2+y^2}$. Without loss of generality, assume that $t'=0$, then
\begin{align}
  t_f=\sqrt{x_{A}^2+y^2}.   \label{eq:bart} 
\end{align}
Thus, using Isaacs' method, we are now able to reduce the solution of the zero-sum differential game of degree to the optimization of a cost/payoff function of one variable. From \eqref{eq:xTt}-\eqref{eq:yDt}, \eqref{eq:phigeo},  and \eqref{eq:bart} we can write the cost (the terminal distance between $A$ and $T$) as $J(y;x_A,x_T,y_T)=\overline{TT_f} \pm \overline{TI}$, where $\overline{TT_f}=\alpha\overline{AI}=\alpha t_f$. Depending on which side of the orthogonal bisector of the segment $AD$ the Target is located we have the two cases:

For case a) where $x_T<0$ solve $\min_y J_1$ where
\begin{align}
     J_1=\alpha\sqrt{x_{A}^2+y^2}+\sqrt{(y_{T}-y)^2+x_{T}^2}.    \label{eq:CostyP}
\end{align}

For case b) where $x_T>0$ solve $\max_y J_2$ where
\begin{align}
     J_2=\alpha\sqrt{x_{A}^2+y^2}-\sqrt{(y-y_{T})^2+x_{T}^2}.    \label{eq:Costy}
\end{align}
 Also, from eq. \eqref{eq:phigeo}  we have that in case a) 
\begin{align}
\left.
	 \begin{array}{l l}
\cos\phi^*=\frac{x_{T}}{\sqrt{x_{T}^2+(y_{T}-y)^2}}, \  \sin\phi^*=\frac{y_{T}-y}{\sqrt{x_{T}^2+(y_{T}-y)^2}}  \label{eq:phia} 
\end{array}  \right.
\end{align}
and in case b)
\begin{align}
\left.
	 \begin{array}{l l}
\cos\phi^*=-\frac{x_{T}}{\sqrt{x_{T}^2+(y-y_{T})^2}}, \  \sin\phi^*\!=\!\frac{y-y_{T}}{\sqrt{x_{T}^2+(y-y_{T})^2}}  \label{eq:phib} 
\end{array}  \right.
\end{align}
The $A$ and $D$ optimal headings are
\begin{align}
  &\cos\chi^*=-\frac{x_{A}}{\sqrt{x_{A}^2+y^2}}, \ \ \ \ \  \sin\chi^*=\frac{y}{\sqrt{x_{A}^2+y^2}}.  \label{eq:chi}   \\
	&  \cos\psi^*=\frac{x_{A}}{\sqrt{x_{A}^2+y^2}}, \ \ \ \ \ \ \sin\psi^*=\frac{y}{\sqrt{x_{A}^2+y^2}}.  \label{eq:psi} 
\end{align}
Let us now use eqs. \eqref{eq:psico} and \eqref{eq:psi} to write the following relationships 
\begin{align}
\left.
	 \begin{array}{l l}
  -\frac{\lambda_{x_D}}{\sqrt{\lambda_{x_D}^2+\lambda_{y_D}^2}}=\frac{x_{A}}{\sqrt{x_{A}^2+y^2}}, \  -\frac{\lambda_{y_D}}{\sqrt{\lambda_{x_D}^2+\lambda_{y_D}^2}}=\frac{y}{\sqrt{x_{A}^2+y^2}}.  \label{eq:Dcostate}  
\end{array}  \right.
\end{align}
Similarly, from eqs. \eqref{eq:chico} and \eqref{eq:chi} we obtain
\begin{align}
  \left.
	 \begin{array}{l l}
	 \frac{\lambda_{x_A}}{\sqrt{\lambda_{x_A}^2+\lambda_{y_A}^2}}=-\frac{x_{A}}{\sqrt{x_{A}^2+y^2}}, \  \frac{\lambda_{y_A}}{\sqrt{\lambda_{x_A}^2+\lambda_{y_A}^2}}=\frac{y}{\sqrt{x_{A}^2+y^2}}.  \label{eq:Acostate}  
\end{array}  \right.  
\end{align}
We have four equations \eqref{eq:Lxad}, \eqref{eq:Lyad}, \eqref{eq:Dcostate}, and \eqref{eq:Acostate} in the four unknowns $\lambda_{x_A}$, $\lambda_{y_A}$, $\lambda_{x_D}$, and $\lambda_{y_D}$. The solution is
\begin{align}
  &\lambda_{x_A}= \frac{1}{2}[x_T(t_f)-x_A(t_f) - \frac{x_{A}}{y}\big(y_T(t_f)-y_A(t_f)\big)]    \nonumber  \\  
  &\lambda_{y_A}= \frac{1}{2}[y_T(t_f)-y_A(t_f) - \frac{y}{x_{A}}\big(x_T(t_f)-x_A(t_f)\big)]    \nonumber  \\  
  &\lambda_{x_D}= \frac{1}{2}[x_T(t_f)-x_A(t_f) + \frac{x_{A}}{y}\big(y_T(t_f)-y_A(t_f)\big)]    \nonumber  \\  
  &\lambda_{y_D}= \frac{1}{2}[y_T(t_f)-y_A(t_f) + \frac{y}{x_{A}}\big(x_T(t_f)-x_A(t_f)\big)]    \nonumber    
\end{align}
By substituting $y_A(t_f)=y$ and $x_A(t_f)=0$ we obtain
\begin{align}
  &\lambda_{x_A}= \frac{1}{2}[x_T(t_f) - \frac{x_{A}}{y}\big(y_T(t_f)-y\big)]    \label{eq:Lxat}  \\  
  &\lambda_{y_A}= \frac{1}{2}[y_T(t_f)-y - \frac{y}{x_{A}}x_T(t_f)]    \label{eq:Lyat}  \\  
  &\lambda_{x_D}= \frac{1}{2}[x_T(t_f) + \frac{x_{A}}{y}\big(y_T(t_f)-y\big)]    \label{eq:Lxdt}  \\  
  &\lambda_{y_D}= \frac{1}{2}[y_T(t_f)-y + \frac{y}{x_{A}}x_T(t_f)]    \label{eq:Lydt}    
\end{align}
which together with eqs. \eqref{eq:Lxtt} and \eqref{eq:Lytt} specify the co-states in terms of the states at time $t_f$. Now, eqs. \eqref{eq:xTt}, \eqref{eq:yTt}, \eqref{eq:bart}, and \eqref{eq:phia} yield
\begin{align}
 a)\  \left\{\left.
	 \begin{array}{l l}
     x_T(t_f)= \Big[1 + \alpha\frac{\sqrt{x_{A}^2+y^2}}{\sqrt{x_{T}^2+(y_{T}-y)^2}}\Big]x_{T} \\
	  y_T(t_f)= \Big[1 + \alpha\frac{\sqrt{x_{A}^2+y^2}}{\sqrt{x_{T}^2+(y_{T}-y)^2}}\Big](y_{T}\!-\!y)\!+\!y 
	\end{array}   \right.  \right.\label{eq:yTopta}  
\end{align}
Similarly, eqs. \eqref{eq:xTt}, \eqref{eq:yTt}, \eqref{eq:bart}, and \eqref{eq:phib} yield
\begin{align}
 b)\  \left\{\left.
	 \begin{array}{l l}
    x_T(t_f)= \Big[1 - \alpha\frac{\sqrt{x_{A}^2+y^2}}{\sqrt{x_{T}^2+(y-y_{T})^2}}\Big]x_{T} \\
	  y_T(t_f)= \Big[1 - \alpha\frac{\sqrt{x_{A}^2+y^2}}{\sqrt{x_{T}^2+(y-y_{T})^2}}\Big](y_{T}\!-\!y)\!+\!y 
	\end{array}   \right.  \right.\label{eq:yToptb}  
\end{align}
Inserting eqs. \eqref{eq:yTopta} and \eqref{eq:yToptb} into eqs. \eqref{eq:Lxtt}, \eqref{eq:Lytt}, and \eqref{eq:Lxat}-\eqref{eq:Lydt} allows us to express the co-states in terms of $y$ (or $t_f$):
\begin{align}
\left.
	 \begin{array}{l l}
   \!\!\!\!\!\!  \lambda_{x_T}= -\Big[1 \pm \alpha\frac{\sqrt{x_{A}^2+y^2}}{\sqrt{x_{T}^2+(y_{T}-y)^2}}\Big]x_{T} \\
 \!\!\!\!\!\! 	  \lambda_{y_T}= -\Big[1 \pm \alpha\frac{\sqrt{x_{A}^2+y^2}}{\sqrt{x_{T}^2+(y_{T}-y)^2}}\Big](y_{T}-y) \\
 \!\!\!\!\!\! 	  \lambda_{x_A}= \frac{1}{2}\Big[1 \pm \alpha\frac{\sqrt{x_{A}^2+y^2}}{\sqrt{x_{T}^2+(y_{T}-y)^2}}\Big](x_{A}+x_{T}-x_{A}\frac{y_{T}}{y}) \\
 \!\!\!\!\!\! 	\lambda_{y_A}= \frac{1}{2}\Big[1 \pm \alpha\frac{\sqrt{x_{A}^2+y^2}}{\sqrt{x_{T}^2+(y_{T}-y)^2}}\Big][y_{T}-(1+\frac{x_{T}}{x_{A}})y] \\
 \!\!\!\!\!\!  \lambda_{x_D}= \frac{1}{2}\Big[1 \pm \alpha\frac{\sqrt{x_{A}^2+y^2}}{\sqrt{x_{T}^2+(y_{T}-y)^2}}\Big](x_{T}-x_{A}+x_{A}\frac{y_{T}}{y}) \\
 \!\!\!\!\!\!  \lambda_{y_D}= \frac{1}{2}\Big[1 \pm \alpha\frac{\sqrt{x_{A}^2+y^2}}{\sqrt{x_{T}^2+(y_{T}-y)^2}}\Big][y_{T}-(1-\frac{x_{T}}{x_{A}})y] 
\end{array}   \right.   \label{eq:costateopta}  
\end{align}
We now used the condition on the Hamiltonian at the terminal time in conjunction with eqs. \eqref{eq:phia}-\eqref{eq:psi} for the optimal controls and eq. \eqref{eq:costateopta} for the co-states. Doing so we obtain the following
\begin{align}
\left.
	 \begin{array}{l l}
-\frac{1}{2}(x_{A}+x_{T}-x_{A}\frac{y_{T}}{y})\frac{x_{A}}{\sqrt{x_{A}^2+y^2}} \\
 + \frac{1}{2}[y_{T}-(1+\frac{x_{T}}{x_{A}})y]\frac{y}{\sqrt{x_{A}^2+y^2}}  \\
  + \frac{1}{2}(x_{T}-x_{A}+x_{A}\frac{y_{T}}{y})\frac{x_{A}}{\sqrt{x_{A}^2+y^2}}  \\
 +\frac{1}{2}[y_{T}-(1-\frac{x_{T}}{x_{A}})y]\frac{y}{\sqrt{x_{A}^2+y^2}} \\
 + \frac{\alpha}{\sqrt{x_{T}^2+(y_{T}-y)^2}}[x_{T}^2+(y_{T}-y)^2]  = 0   \\
\Rightarrow  \ \  -x_{A}^2 + x_{A}^2\frac{y_{T}}{y} + y_{T}y -y^2  \\
\qquad  +\alpha \frac{\sqrt{x_{A}^2+y^2}}{\sqrt{x_{T}^2+(y_{T}-y)^2}}  [x_{T}^2+(y_{T}-y)^2] = 0    \\
\Rightarrow  \ \ (x_{A}^2 + y^2 - x_{A}^2\frac{y_{T}}{y} - y_{T}y)^2 \\
\qquad  = \alpha^2 (x_{A}^2+y^2) [x_{T}^2+(y_{T}-y)^2]  \\
\Rightarrow  \ \ (1 - \frac{y_{T}}{y})^2(x_{A}^2+y^2) = \alpha^2 [x_{T}^2+(y_{T}-y)^2] 
\end{array}   \right.   \nonumber
\end{align}
Finally, grouping terms in $y$ we obtain \eqref{eq:Quartic}, that is, the optimal aimpoint on the orthogonal bisector of $\overline{AD}$ specified by $y^*$ is a real solution of the quartic equation \eqref{eq:Quartic}.  

For the derivation of the Target strategy when $x_T=0$, please refer to  \cite{Pachter19}.  $\square$  

The verification of the saddle-point strategies of the ATTDG in the escape region was recently provided in \cite{garcia2019}.
Preliminary results and a candidate solution for the ATDDG in the capture region have been presented in \cite{garcia18Optimal}; 
verification of these strategies is a topic of current research.

\section{Concluding Remarks}
Pursuit-evasion differential games, initially introduced by Rufus Isaacs, have been studied in great detail over the past half-century. Since its inception, many ideas and techniques have been developed to gain a better understanding of pursuit-evasion differential games. In his seminal work, Isaacs motivated and described the idea of posing problems in a game-theoretic framework; the dynamics are mathematically described by differential equations so he called this approach ``Differential Games''. Using differential games, problems of pursuit-evasion have been properly described and solved. The contributions of Rufus Isaacs, Richard Bellman, John Breakwell and his students, and Lev S. Pontryagin were highlighted in this paper.
A description of different categories of pursuit-evasion differential games and significant contributions were provided with examples. Classic games such as the ``Homicidal Chauffeur'' Differential Game and ``Two-Cars'' Differential Game were briefly discussed. Games involving N-pursuers-1-evader, 1-pursuer-M-evaders, and N-pursuer-M-evaders were addressed. In games with multiple players, cooperation between members of the same team was emphasized. Not only cooperation improves team performance but it was shown that cooperative behaviors are necessary to synthesize saddle-point equilibrium strategies. Two pursuit-evasion problems involving multiple players were presented and they were formulated as differential games. The solutions provided in this paper highlighted the cooperative aspect and the methods available by applying differential game theory.

\bibliographystyle{IEEEtran}
\bibliography{references}

\end{document}